\documentclass[11pt]{article}

\usepackage{amsmath}
\usepackage{amsthm}
\usepackage{amsfonts}
\usepackage{fancybox}
\usepackage{epsfig}
\usepackage{epsf}
\usepackage{amssymb}
\usepackage{epic,eepic}
\usepackage{latexsym,bm}
\usepackage{graphicx}
\usepackage{float}
\usepackage{color}
\usepackage[top=1in, bottom=1in, left=1in, right=1in, dvips,letterpaper]{geometry}
\usepackage{algorithm}
\usepackage{algorithmic}
\usepackage{cite}
\usepackage{boxedminipage}
\usepackage{appendix}
\usepackage{sectsty}

\newcommand\argmin{\mathop{\textrm{argmin}}}
\newcommand\crule[1][5cm]{%
  \par
  \nointerlineskip
  \centerline{\hbox to #1{\hrulefill}}%
  \nointerlineskip}

\numberwithin{equation}{section}
\numberwithin{algorithm}{section}
\newtheorem{thm}{{\sc Theorem}}[section]
\newtheorem{lem}{{\sc Lemma}}[section]

\newtheorem{rem}{Remark}[section]

\newtheorem{defi}{{\sc Definition}}[section]

\newcommand{\R}{\mathbb{R}}
\newcommand{\E}{\mathbb{E}}
\newcommand{\SFO}{\mathcal{SFO}}
\newcommand{\SZO}{\mathcal{SZO}}

\newcommand{\be}{\begin{equation}}
\newcommand{\ee}{\end{equation}}
\newcommand{\bee}{\begin{equation*}}
\newcommand{\eee}{\end{equation*}}
\newcommand{\bea}{\begin{eqnarray}}
\newcommand{\eea}{\end{eqnarray}}
\newcommand{\beaa}{\begin{eqnarray*}}
\newcommand{\eeaa}{\end{eqnarray*}}

\newcommand{\etal}{{\it et al.\ }}
\newcommand{\st}{\textnormal{s.t.}}

\title{Penalty Methods with Stochastic Approximation for Stochastic Nonlinear Programming}

\author{Xiao Wang
\thanks{School of Mathematical Sciences, University of Chinese Academy of Sciences; Key Laboratory of Big Data Mining and Knowledge Management, Chinese Academy of Sciences, China. Email: wangxiao@ucas.ac.cn. Research of this author was supported in part by Postdoc Grant 119103S175, UCAS President Grant Y35101AY00 and NSFC Grant 11301505.}
\and Shiqian Ma
\thanks{Department of Systems Engineering and Engineering Management, The Chinese University of Hong Kong, Shatin, N. T., Hong Kong. Email: sqma@se.cuhk.edu.hk. Research of this author was supported in part by a Direct Grant of the Chinese University of Hong
Kong (Project ID: 4055016) and the Hong Kong Research Grants Council
General Research Fund Early Career Scheme (Project ID: CUHK 439513).}
\and Ya-xiang Yuan
\thanks{State Key Laboratory of Scientific and Engineering Computing, Academy of Mathematics and Systems Science, Chinese Academy of Sciences,
China. Email: yyx@lsec.cc.ac.cn. Research of this author was supported in part by NSFC Grants 11331012 and 11321061.}}

\date{May 18, 2016}

\begin{document}

\maketitle

\begin{abstract}
In this paper, we propose a class of penalty methods with stochastic approximation for solving stochastic nonlinear programming problems. We assume that only noisy gradients or function values of the objective function are available via calls to a stochastic first-order or zeroth-order oracle. In each iteration of the proposed methods, we minimize an exact penalty function which is nonsmooth and nonconvex with only stochastic first-order or zeroth-order information available. Stochastic approximation algorithms are presented for solving this particular subproblem. The worst-case complexity of calls to the stochastic first-order (or zeroth-order) oracle for the proposed penalty methods for obtaining an $\epsilon$-stochastic critical point is analyzed.

\vspace{0.8cm}


\vspace{0.5cm}

\end{abstract}

\pagestyle{myheadings}
\thispagestyle{plain}


\section{Introduction}

In this paper, we consider the following stochastic nonlinear programming (SNLP) problem:
\begin{equation}\label{orig-prob}
\begin{split}
\min_{x\in\R^n}\quad & f(x)\\
\st \quad & c(x):=(c_1(x),\ldots,c_q(x))^T=0,
\end{split}
\end{equation}
where both $f:\R^n \to \R$ and $c:\R^n\to \R^q$ are continuously differentiable but possibly nonconvex.
We assume that the function values and gradients of $c_i(x)$, $i=1,\ldots,q$, can be obtained exactly. However, we assume that only the noisy function values or gradients of $f$ are available. Specifically, the noisy gradients (resp. function values) of $f$ are obtained via subsequent calls to a {\it stochastic first-order oracle} ($\SFO$) (resp. {\it stochastic zeroth-order oracle} ($\SZO$)). The problem \eqref{orig-prob} arises in many applications, such as machine learning \cite{mbps09}, simulation-based optimization \cite{f02}, mixed logit modeling problems in economics and transportation \cite{bbt00,bct06,hg03}. Besides, many two-stage stochastic programming problems can be formulated as \eqref{orig-prob} (see, e.g., \cite{bl11}). Many problems in these fields have the following objective functions:
\[
f(x) = \int_{\Xi} F(x,\xi)dP(\xi) \quad \mbox{or} \quad f(x)=\E_\xi[F(x,\xi)],
\]
where $\xi$ denotes the random variable whose distribution $P$ is supported on $\Xi$ and $\E_{\xi}[\cdot]$ means that the expectation is taken with respect to $\xi$. Due to the fact that the integral is difficult to evaluate, or function $F(\cdot,\xi)$ is not given explicitly, the function values and gradients of $f$ are not easily obtainable and only noisy information of $f$ is available.

Stochastic programming has been studied for several decades. Robbins and Monro \cite{rm51} proposed a stochastic approximation (SA) algorithm for solving convex stochastic programming problems. Various methods on SA have been proposed after \cite{rm51}, such as \cite{C54,e83,g78,rs86,s58} and so on. By incorporating the averaging technique, Polyak \cite{p90} and Polyak and Juditsky \cite{pj92} suggested SA methods with longer stepsizes and the {\it asymptotically optimal} rate of convergence is exhibited.
Interested readers are referred to \cite{bl11,sdr09} for more details on stochastic programming.
Recently, following the development of the complexity theory in convex optimization  \cite{ny83}, the convergence and complexity properties of SA methods were explored. Nemirovski \etal \cite{njls09} proposed a mirror descent SA method for the nonsmooth convex stochastic programming problem $x^*:=\argmin \{ f(x)\mid x\in X\}$ and showed that the algorithm returns $\bar{x}\in X$ with $\E[f(\bar{x})-f(x^*)]\le\epsilon$ in $O(\epsilon^{-2})$ iterations, where $X$ is the constraint set and $\E[y]$ denotes the expectation of random variable $y$. Nemirovski and Rubinstein \cite{NR} proposed an efficient SA method for convex-concave stochastic saddle point problem in the form of $\min_{x\in X}\max_{y\in Y}\phi(x,y)$. It is assumed that both $X$ and $Y$ are convex sets and $\phi$ is convex in $x\in X$ and concave in $y\in Y$. Under certain assumptions, they showed that the proposed method returns $(\bar{x},\bar{y})\in X\times Y$ with $\E[\max_{y\in Y}\phi(\bar{x},y)-\min_{x\in X}\phi(x,\bar{y})]\le\epsilon$ in $O(\epsilon^{-2})$ iterations.
Recently, Wang and Bertsekas \cite{WB} proposed an SA method with constraint projection for nonsmooth convex optimization, whose constraint set is the intersection of a finite number of convex sets. Other relevant works on the complexity analysis of SA algorithms for convex optimization include \cite{gl12,gl131,jrt08,ksm01,l12,lns12}.

SA algorithms for nonconvex stochastic programming and their complexity analysis, however, have not been investigated thoroughly yet.
In \cite{gl13}, Ghadimi and Lan proposed an SA method for the nonconvex stochastic optimization problem $\min\{f(x)\mid x\in \R^n\}$. Their algorithm returns $\bar{x}$ with $\E[\|\nabla f(\bar{x})\|^2]\le \epsilon$ after at most $O(\epsilon^{-2})$ iterations. In \cite{glz13}, Ghadimi \etal studied the following nonconvex composite stochastic programming problem
\be\label{glz-prob}
\min_{x\in X}\quad f(x) + \ell(x),
\ee
where $X\subseteq\R^n$ is a closed convex set, $f$ is nonconvex and $\ell$ is a simple convex function with certain special structure. They proposed a proximal-gradient like SA method for solving \eqref{glz-prob} and analyzed its complexity.
Dang and Lan \cite{dl13} studied several stochastic block mirror descent methods for large-scale nonsmooth and stochastic optimization by combining the block-coordinate decomposition and an incremental block average scheme. In \cite{gl132}, Ghadimi and Lan generalized Nesterov's accelerated gradient method \cite{Nesterov-1983} to solve the stochastic composite optimization problem \eqref{glz-prob} with $X:=\R^n$. However, to the best of our knowledge, there has not been any SA method proposed for solving SNLP \eqref{orig-prob} with nonconvex objective functions and nonconvex constraints. In this paper, we will focus on studying such methods and analyzing their complexity properties.

When the exact gradient of $f$ in \eqref{orig-prob} is available, a classical way to solve \eqref{orig-prob} is using penalty methods.
In a typical iteration of a penalty method for solving \eqref{orig-prob}, an associated penalty function is minimized for a fixed penalty parameter. The penalty parameter is then adjusted for the next iteration. For example, the exact penalty function $\Phi_\rho(x)  = f(x) + \rho\|c(x)\|_2$ is widely used in penalty methods (see, e.g., \cite{cgt11}). Note that $\Phi_{\rho}$ is the summation of a differentiable term and a nonsmooth term, and the nonsmooth term itself is the composition of the convex nonsmooth function $\rho\|\cdot\|_2$ and a nonconvex differentiable function $c(x)$. In \cite{cgt11}, an exact penalty algorithm is proposed for solving \eqref{orig-prob} which minimizes $\Phi_\rho(x)$ in each iteration with varying $\rho$, and its function-evaluation worst-case complexity is analyzed. We refer the interested readers to \cite{nw06} for more details on penalty methods.

Motivated by the work in \cite{cgt11}, we shall propose a class of penalty methods with stochastic approximation in this paper for solving SNLP \eqref{orig-prob}. In our methods, we minimize a penalty function $f(x)+\rho\|c(x)\|_2$ in each iteration with varying $\rho$. Note that the difference is that now we only have access to inexact information to $f$ through $\SFO$ or $\SZO$ calls. We shall show that our proposed methods can return an $\epsilon$-stochastic critical point (will be defined later) of \eqref{orig-prob}, and {analyze the worst-case complexity of $\SFO$ (or $\SZO$) calls to obtain such a solution}.

{\bf Contributions.} Our contributions in this paper lie in the following folds. 
First, we propose a penalty method with stochastic first-order information for solving \eqref{orig-prob}. In each iteration of this algorithm, we solve a nonconvex stochastic composite optimization problem as a subproblem. {An SA} algorithm for solving this subproblem is also given. The $\SFO$-calls worst-case complexity of this penalty method to obtain an $\epsilon$-stochastic critical point is analyzed. Second, for problem \eqref{orig-prob} with only stochastic zeroth-order information (i.e., noisy function values) available, we also present a penalty method for solving them and analyze their $\SZO$-calls worst-case complexity.

{\bf Notation.} We adopt the following notation throughout the paper. $\nabla f(x)$ denotes the gradient of $f$ and ${J}(x):=\nabla c(x) = (\nabla c_1(x),\ldots,\nabla c_q(x))^T$ denotes the Jacobian matrix of $c$. The subscript $_k$ refers to the iteration number in an algorithm, e.g., $x_k$ is the $k$-th $x$ iterate. {$x^Ty$} denotes the Euclidean inner product of vectors $x$ and $y$ in $\R^n$. Without specification, $\|\cdot\|$ represents the Euclidean norm $\|\cdot\|_2$ in $\R^n$.

{\bf Organization.} The rest of this paper is organized as follows. In Section \ref{sec:sco-1st}, we propose an SA algorithm with stochastic first-order information for solving a nonconvex stochastic composite optimization problem \eqref{nonsmooth-uncons}, which is the subproblem in our penalty methods for solving \eqref{orig-prob}. In Section \ref{sec:penalty-1st}, we propose a penalty method with stochastic first-order information for solving the {SNLP} problem \eqref{orig-prob} and analyze its $\SFO$-calls worst-case complexity to obtain an $\epsilon$-stochastic critical point. In Section \ref{sec:penalty-0th}, we present a penalty method with SA for solving \eqref{orig-prob} using only stochastic zeroth-order information of $f$ and analyze its $\SZO$-calls worst-case complexity. Finally, we draw some conclusions in Section \ref{sec:conclusions}.


\section{A stochastic first-order approximation method for a nonconvex stochastic composite optimization}\label{sec:sco-1st}

Before we present the penalty methods for solving SNLP \eqref{orig-prob}, we consider the following nonconvex stochastic composite optimization (NSCO) problem in this section, which is in fact the subproblem in our penalty methods for solving \eqref{orig-prob}:
\be \label{nonsmooth-uncons}
\min_{x\in\R^n}\quad \Phi_h(x) := f(x) + h(c(x)),
\ee
where $f$ and $c$ are both continuously differentiable and possibly nonconvex, and $h$ is a nonsmooth convex function. We assume that both the exact zeroth-order and first-order information (function value and Jacobian matrix) of $c$ is available, but only noisy gradient information of $f$ is available via $\SFO$ calls. Namely, for the input $x$, $\SFO$ will output a stochastic gradient $G(x,\xi)$ of $f$, where $\xi$ is a random variable whose distribution is supported on $\Xi\subseteq \R^d$ (note that $\Xi$ does not depend on $x$).

NSCO \eqref{nonsmooth-uncons} is quite different from \eqref{glz-prob} considered by Ghadimi \etal in \cite{glz13}. In \eqref{glz-prob}, the second term in the objective function must be convex. However, we allow $c(x)$ to be nonconvex which implies that the second term $h(c(x))$ in \eqref{nonsmooth-uncons} is nonconvex. For solving \eqref{nonsmooth-uncons} under deterministic settings, i.e., when exact zeroth-order and first-order information of $f$ is available, there have been some relevant works. Cartis \etal \cite{cgt11} proposed a trust region approach and a quadratic regularization approach for solving \eqref{nonsmooth-uncons}, and explored their function-evaluation worst-case complexity. Both methods need to take at most $O(\epsilon^{-2})$ function-evaluations to reduce a first-order criticality measure below $\epsilon$. Garmanjani and Vicente \cite{gv12} proposed a smoothing direct-search method for nonsmooth nonconvex but Lipschitz continuous unconstrained optimization. They showed that the method takes at most $O(\epsilon^{-3}\log \epsilon^{-1})$ function-evaluations to reduce both the smoothing parameter and the first-order criticality of the smoothing function below $\epsilon$. Bian and Chen \cite{bc13} studied the worst-case complexity of a smoothing quadratic regularization method for a class of nonconvex, nonsmooth and non-Lipschitzian unconstrained optimization problems. Specifically, by assuming $h(c(x)):=\sum_{i=1}^n\phi(|x_i|^p)$ in \eqref{nonsmooth-uncons}, where $0<p\le 1$ and $\phi$ is some continuously differentiable function, it was shown in \cite{bc13} that the function-evaluation worst-case complexity to reach an $\epsilon$ scaled critical point is $O(\epsilon^{-2})$. However, to the best of our knowledge, there has not been any work studying NSCO \eqref{nonsmooth-uncons}.

The following assumptions are made throughout this paper.
\begin{description}
\item[]\textbf{AS.1}\quad $f,c_i\in\mathcal{C}^1(\R^n)$ \footnote{$f\in\mathcal{C}^{1}(\R^n)$ means that $f: \R^n\to \R$ is continuously differentiable.}, $i=1, \ldots, q$. $f(x)$ is lower bounded by {a real number} $f^{low}$ for any $x\in\R^n$. $\nabla f$ and $J$ are Lipschitz continuous with Lipschitz constants $L_g$ and $L_J$ respectively.
\item[]\textbf{AS.2}\quad $h$ is convex and Lipschitz continuous with Lipschitz constant $L_h$.
\item[]{\bf AS.3}\quad $\Phi_h(x)$ is lower bounded by {a real number} $\Phi_h^{low}$ for all $x\in\R^n$.
\item[]\textbf{AS.4}\quad For any $k$, we have
\begin{align*}
  a)\quad & \E\left[G(x_k,\xi_k)\right]=\nabla f(x_k),\\
  b)\quad & \E\left[\|G(x_k,\xi_k)-\nabla f(x_k)\|^2\right]\le \sigma^2,
\end{align*}
where $\sigma>0$.
\end{description}

We now describe our {SA} algorithm for solving NSCO \eqref{nonsmooth-uncons} in Algorithm \ref{alg-uncons}. For ease of presentation, we denote
\begin{equation}\label{def-psi-gamma}
\psi_\gamma (x,g,u):= {g^T(u-x)} + h(c(x)+{J}(x)(u-x)) + \frac{1}{2\gamma}\|u-x\|^2.
\end{equation}

\begin{algorithm}[H]
\caption{\quad \bf Stochastic approximation algorithm for NSCO \eqref{nonsmooth-uncons}}
\label{alg-uncons}
\begin{algorithmic}[1]
\REQUIRE Given $x_1\in\R^n$, maximum iteration number {$N_{in}$}, stepsizes $\{\gamma_k\}$ with $\gamma_k>0$, $k\ge1$, the batch sizes $\{m_k\}$ with $m_k>0$, $k\ge1$. Let $R$ be a random variable following probability distribution $P_R$ which is supported on $\{1,\ldots,N_{in}\}$.
\ENSURE {$x_R$}.
\FOR{$k=1,2,\ldots, R-1$}
\STATE Call $\SFO$ $m_k$ times to obtain $G(x_k,\xi_{k,i})$, $i=1,\ldots,m_k$, then set
    \begin{equation*}
    G_k = \frac{1}{m_k}\sum_{i=1}^{m_k}G(x_k,\xi_{k,i}).
    \end{equation*}
\STATE Compute
\begin{equation}\label{x_k+1}
x_{k+1} = \argmin_{u\in\R^n} \psi_{\gamma_k} (x_k,G_k,u).
\end{equation}
\ENDFOR
\end{algorithmic}
\end{algorithm}

The most significant difference between our strategy to update iterates in \eqref{x_k+1} and the one in \cite{glz13} is the way that we deal with the structured nonsmooth term $h(c(x))$. Since it is the composition of the nonsmooth convex function $h$ and the nonconvex differentiable function $c$, we apply the first-order approximation of $c$ in \eqref{x_k+1}. Due to the convexity of $h$, $\psi_\gamma$ is strongly convex with respect to $u$. Hence, $x_{k+1}$ is well-defined in \eqref{x_k+1}.

Let us define
\begin{equation}\label{proj-gra}
P_\gamma (x,g): = \frac{1}{\gamma}(x-x^+),
\end{equation}
where $x^+$ is defined as
\begin{equation}\label{x_+}
x^+ = \argmin_{u\in\R^n} \ \psi_\gamma (x,g,u).
\end{equation}
From the optimality conditions for \eqref{x_+}, it follows that there exists
$
p\in\partial h(c(x)+{J}(x)({x^+}-x))
$ such that $P_\gamma(x,g)=g + {J}(x)^Tp$. Thus,
if $P_\gamma(x,\nabla f(x))=0$, then $x$ is a first-order critical point of \eqref{nonsmooth-uncons}. Therefore, $\|P_\gamma(x,\nabla f(x))\|$ can be adopted as the criticality measure for \eqref{nonsmooth-uncons}. In addition, we denote the generalized gradients
\be \label{g-rand_g}
\tilde{g}_k := P_{\gamma_k}(x_k, \nabla f(x_k))\quad \mbox{and}\quad \tilde{g}_k^r := P_{\gamma_k}(x_k,G_k).
\ee
The following results give estimates to $\E[\|\tilde{g}_R\|^2]$ and $\E[\|\tilde{g}_R^r\|^2]$.

As the analysis in this section essentially follows from \cite{glz13}, for simplicity we only state the results here and their proofs are given in Appendix \ref{appA}. The first theorem provides an upper bound for the expectation of the generalized gradient at $x_R$, the output of Algorithm \ref{alg-uncons}.
\begin{thm}\label{thm2.2}
Let {\bf AS.1-4} hold. We assume that the stepsizes $\{\gamma_k\}$ in Algorithm \ref{alg-uncons} are chosen such that $0<\gamma_k\le 2/L$ with $\gamma_k<2/L$ for at least one $k$, where $L:=L_g+L_hL_J$. Moreover, suppose that the probability mass function $P_R$ is chosen such that for any $k=1,\ldots,N_{in}$,
\begin{equation}\label{prob_R}
P_R(k):=\mathrm{Prob}\{R=k\}=\frac{ \gamma_k - L\gamma_k^2/2}{\sum_{k=1}^{{N_{in}}}( \gamma_k-L\gamma_k^2/2)}.
\end{equation}
Then for any ${N_{in}}\ge1$, we have
\be\label{proj-gra-rand-bound}
\E[\|\tilde{g}_R^r\|^2]\le \frac{D_{\Phi_h} +  \sigma^2 \sum_{k=1}^{N_{in}}(\gamma_k/m_k)}{\sum_{k=1}^{N_{in}}( \gamma_k - L\gamma_k^2/2)},
\ee
where the expectation is taken with respect to $R$ and $\xi_{[N_{in}]}:=(\xi_1,\ldots,\xi_{N_{in}})$ with $\xi_k:=(\xi_{k,1},\ldots,\xi_{k,m_k})$ in Algorithm \ref{alg-uncons}, and the real number $D_{\Phi_h}$ is defined as
\be \label{D_Phih}
D_{\Phi_h}=\Phi_h(x_1)-\Phi_h^{low}.
\ee
\end{thm}

By specializing the settings of Algorithm \ref{alg-uncons}, we obtain the following complexity result.

\begin{thm}\label{cor2.5}
Let {\bf AS.1-4} hold. Suppose that in Algorithm \ref{alg-uncons}, $\gamma_k=1/L$ where $L:=L_g + L_h L_J$ and the probability mass function is chosen as in \eqref{prob_R}. For any given $\epsilon>0$, we assume that the total number of $\SFO$-calls $\bar{N}$ in Algorithm \ref{alg-uncons} satisfies
\be\label{bound-N}
\bar{N} \ge \max\left\{\frac{\left(D_{\Phi_h}C_2 + L C_3\right)^2}{\epsilon^2} + \frac{32LD_{\Phi_h}}{\epsilon}, \frac{C_1}{L^2}\right\},
\ee
where
\be \label{C_1-C_2-C_3}
C_1=\sigma^2/\tilde{D},\quad C_2=8\sigma/\sqrt{\tilde{D}}\quad \mbox{and}\quad C_3=6\sigma\sqrt{\tilde{D}}
\ee
with some problem-independent positive constant $\tilde{D}$.
We further assume that the batch size $m_k$, $k=1,\ldots,N_{in}$, satisfies
\be \label{batch-size-m}
m_k = m := \left\lceil \min \left\{ \bar{N}, \max\left\{1, \frac{\sigma}{L}\sqrt{\frac{\bar{N}}{\tilde{D}}}\right\}\right\}\right\rceil,
\ee
Then we have
\be
\E[\|\tilde{g}_R\|^2]\le\epsilon \quad \mbox{and}\quad \E[\|\tilde{g}_R^r\|^2]\le\epsilon,
\ee
where the expectations are taken with respect to $R$ and $\xi_{[N_{in}]}$. Thus, it follows that the number of $\SFO$-calls required by Algorithm \ref{alg-uncons} to achieve $\E[\|\tilde{g}_R\|^2]\le\epsilon$ and $\E[\|\tilde{g}_R^r\|^2]\le\epsilon$ is in the order of $O(\epsilon^{-2})$.
\end{thm}

\begin{rem}
{Theorems \ref{thm2.2}-\ref{cor2.5} are similar to the theoretical results obtained in \cite{glz13}. The only difference is that we allow the nonsmooth term of the objective to be nonconvex, while the results in \cite{glz13} require the nonsmooth term to be convex.}
\end{rem}

\begin{rem}
{Instead of choosing a random iterate as the output, we can use a deterministic termination condition, i.e., choosing the iterate  $\hat{x}$ that has the smallest norm of the exact gradient among all iterates as the output of the algorithm. Following the analysis in Theorems \ref{thm2.2}-\ref{cor2.5}, we can obtain a similar bound on the expectation of the squared norm of the gradient at $\hat{x}$ and obtain the same complexity result $O(\epsilon^{-2})$. However, this deterministic termination condition requires to compute the exact gradients at all iterates, which is impractical for stochastic programming.}
\end{rem}

%



\section{A penalty method with stochastic first-order approximation for SNLP \eqref{orig-prob}}\label{sec:penalty-1st}

We now return to the {SNLP} problem \eqref{orig-prob}, in which only stochastic gradient information of $f$ is available via $\SFO$-calls. In this section, we shall propose a penalty method with stochastic {first-order} approximation for solving \eqref{orig-prob} and study its $\SFO$-calls worst-case complexity.

In deterministic settings, one would expect to find the KKT point of \eqref{orig-prob}, which is defined as follows (see \cite{nw06} for reference).
\begin{defi}\label{KKT-def}
$x^*$ is called a KKT point of \eqref{orig-prob}, if there exists $\lambda^*\in\R^q$ such that
\begin{equation*}
\nabla f(x^*) + J(x^*)^T\lambda^* = 0, \mbox{ and }\quad c(x^*) = 0.
\end{equation*}
\end{defi}

When solving nonlinear programming {problems}, however, it is possible that one algorithm fails to output a feasible point. For example, the constraints $c(x)=0$ may not be realized for any $x\in\R^n$. In this case, the best one can hope is to find $x$ such that $\|c(x)\|$ is minimized, or in other words, the constraint violation $\|c(x)\|$ could not be improved any more in a neighborhood of $x$. Therefore, Cartis, Gould and Toint \cite{cgt11} introduced the following definition of $\epsilon$-approximate critical point of \eqref{orig-prob}.
\begin{defi}\label{def-cartis-epsilon-optimal}
$x$ is called an $\epsilon$-approximate critical point of \eqref{orig-prob}, if there exists $\lambda\in\R^q$ such that the following two inequalities hold:
\begin{equation*}\label{def-cartis-epsilon-optimal-two-conditions}
\|\nabla f(x) + {J}(x)^T\lambda\|\le \epsilon, \mbox{ and }\quad \theta(x)\le \epsilon,
\end{equation*}
where $\theta(x)$ is defined as
\be \label{theta}
\theta(x) = \|c(x)\| - \min_{\|s\|\le 1} \|c(x) + {J}(x)s\|.
\ee
\end{defi}
Note that $\bar{x}$ is a critical point of the problem $\{\min \ \|c(x)\|\}$, if $\theta(\bar{x})=0$ (see e.g. \cite{cgt11,y85}).

In stochastic settings, any specific algorithm for solving \eqref{orig-prob} is a random process and the output is a random variable. We thus modify Definition \ref{def-cartis-epsilon-optimal} and define the {\it $\epsilon$-stochastic critical point} of \eqref{orig-prob} as follows.
\begin{defi} \label{def3.1}
Let $\epsilon$ be any given positive constant and $x\in\R^n$ be output of a random process. $x$ is called an $\epsilon$-stochastic critical point of \eqref{orig-prob}, if there exists $\lambda\in\R^q$ such that
\begin{align}
&\E[\|\nabla f(x) + {J}(x)^T\lambda\|^2] \le \epsilon, \label{opt-cond-1}\\
&\E[\theta(x)] \le \sqrt{\epsilon}. \label{opt-cond-2}
\end{align}
\end{defi}

{
We now make a few remarks regarding to this definition. In the deterministic setting, \eqref{opt-cond-1} and \eqref{opt-cond-2} reduce respectively to $\|\nabla f(x) + {J}(x)^T\lambda\|\le \sqrt{\epsilon}$ and $\theta(x)\le \sqrt{\epsilon}$, which are both worse than the conditions in Definition \ref{def-cartis-epsilon-optimal}. In \eqref{opt-cond-1} we use $\E[\|\nabla f(x) + {J}(x)^T\lambda\|^2]$ instead of $\E[\|\nabla f(x) + {J}(x)^T\lambda\|]$, because for the subproblem NSCO \eqref{nonsmooth-uncons} we are only able to analyze the former term. It is worth noting that by Jensen's inequality, we have $\|\E[\nabla f(x) + {J}(x)^T\lambda]\|^2 \leq \E[\|\nabla f(x) + {J}(x)^T\lambda\|^2]$, and are able to bound $\|\E[\nabla f(x) + {J}(x)^T\lambda]\|$. However, our analysis is directly for $\E[\|\nabla f(x) + {J}(x)^T\lambda\|^2]$, and replacing it by $\|\E[\nabla f(x) + {J}(x)^T\lambda]\|$ in Definition \ref{def3.1} will loosen the bound. 
Admittedly, the bounds in Definition \ref{def3.1} are loose compared with the ones in Definition \ref{def-cartis-epsilon-optimal}. However, note that Definition \ref{def3.1} is for SNLP \eqref{orig-prob} in the stochastic setting, and that is the price we need to pay when we define the $\epsilon$-stochastic critical point.}

We now give our penalty method with stochastic first-order approximation for solving SNLP \eqref{orig-prob}. Similar as the deterministic penalty method in \cite{cgt11}, we minimize, at each iteration, the following penalty function with varying penalty parameter $\rho$:
\be \label{exact-penalty}
\min_{x\in\R^n}\quad \Phi_\rho(x) = f(x) + \rho\|c(x)\|.
\ee
Notice that \eqref{exact-penalty} is a special case of NSCO \eqref{nonsmooth-uncons} with $h(\cdot):=\rho\|\cdot\|$. Hence, $h$ is convex and Lipschitz continuous with Lipschitz constant $L_h=\rho$. {\bf AS.2} thus holds naturally. Moreover, if {\bf AS.1} is assumed to be true, then for any $\rho>0$, there exists $\Phi_\rho^{low} \ge f^{low}$ such that $\Phi_\rho(x)\ge\Phi_\rho^{low}$ for all $x\in\R^n$. Therefore, {\bf AS.3} holds as well with $h(\cdot):=\rho\|\cdot\|$ and $\Phi_h^{low}:=\Phi_\rho^{low}$.
Our penalty method for solving \eqref{orig-prob} is described in Algorithm \ref{alg-cons}. 

\begin{algorithm}[H]
\caption{\quad \bf Penalty method with stochastic first-order approximation for \eqref{orig-prob}}
\label{alg-cons}
\begin{algorithmic}[1]
\REQUIRE Given $N$ as the maximum iteration number, tolerance $\epsilon\in(0,1)$, steering parameter $\xi\in(0,1)$, initial iterate $x_1\in\R^n$, $G_1\in\R^n$, penalty parameter $\rho_0 \ge1$, minimal increase factor $\tau>0$. Set $k:=1$.
\ENSURE {$x_N$}.
\FOR{$k=1,2,\ldots, N-1$}
\STATE Step (a): Find $\rho:=\rho_k\ge\rho_{k-1}+\tau$ satisfying
\be \label{phi-theta}
\phi_\rho(x_k) \ge \rho \xi\theta(x_k),
\ee
where $\theta(x)$ is defined in \eqref{theta} and
\be\label{phi}
\phi_\rho(x_k) = \rho \|c(x_k)\| - \min_{\|s\|\le1}\left\{{G_k^Ts} + \rho\|c(x_k) + {J}(x_k)s\|\right\}.
\ee
\STATE Step (b): Apply Algorithm \ref{alg-uncons} with initial iterate $x_{k,1}:=x_k$ to solve the NSCO subproblem \eqref{exact-penalty} with $\rho:=\rho_k$ and {using $\bar{N}_\rho$ $\SFO$-calls}, returning $x_{k+1}:=x_{k,R_k}$ and $G_{k+1}:=G_{k,R_k}$, such that
\be \label{epsi-g}
\E[\|\tilde{g}_{k+1}^r\|^2] \le \epsilon,
\ee
where $\tilde{g}_k^r$ is defined in \eqref{g-rand_g}, $x_{k,R_k}$ denotes the $R_k$-th iterate generated by Algorithm \ref{alg-uncons} when solving the $k$-th subproblem, and the expectation is taken with respect to the random variables generated when calling Algorithm \ref{alg-uncons}.
\ENDFOR
\end{algorithmic}
\end{algorithm}

Note that Algorithm \ref{alg-cons} provides a unified framework of penalty methods for SNLP \eqref{orig-prob}, and any algorithm for solving NSCO in Step (b) can be incorporated into Algorithm \ref{alg-cons}.

\begin{rem}\label{rem3.1}
We now remark that Step (a) in Algorithm \ref{alg-cons} is well-defined, i.e., \eqref{phi-theta} can be satisfied for sufficiently large penalty parameter $\rho$. This fact can be seen from the following argument:
\begin{align*}
\phi_\rho(x_k) & =  \rho \|c(x_k)\| - \min_{\|s\|\le1}\left\{{G_k^Ts} + \rho\|c(x_k) + {J}(x_k)s\|\right\} \notag \\
               & \ge \rho\|c(x_k)\| - \min_{\|s\|\le1} \left\{\|G_k\| + \rho\|c(x_k) + {J}(x_k)s\|\right\} \notag \\
               & =  -\|G_k\|  + \rho\left\{\|c(x_k)\| - \min_{\|s\|\le 1}\|c(x_k) + {J}(x_k)s\|\right\} \notag \\
               & =  -\|G_k\| + \rho \theta(x_k). 
\end{align*}
This indicates that \eqref{phi-theta} holds when
\be \label{rho-low-bound}
\rho \ge \frac{\|G_k\|}{(1 - \xi)\theta(x_k)}.
\ee
Once the algorithm enters Step (a), both $x_k$ and $G_k$ are fixed, so we can achieve \eqref{rho-low-bound} by increasing $\rho$.
\end{rem}

\begin{rem}
Although motivated by the exact penalty-function algorithm proposed in \cite{cgt11} for solving nonlinear programming in the deterministic setting, our Algorithm \ref{alg-cons}, as {an SA} method, is significantly different from the algorithm in \cite{cgt11} in the following folds.
\begin{itemize}
\item[(i)] Different subproblem solver is used in Algorithm \ref{alg-cons}. In \cite{cgt11}, each composite optimization subproblem is solved by a trust region algorithm {or a quadratic-regularization algorithm}. For stochastic programming, however, since exact objective gradient is not available, exact gradient-based algorithms do not work any more. So we adopt a stochastic approximation algorithm to solve NSCO subproblems in Algorithm \ref{alg-cons}. This will yield quite different subproblem termination criterion.
\item[(ii)] Different termination condition for the subproblem is used in Algorithm \ref{alg-cons}. When subproblems in \cite{cgt11} are solved, an extra condition $\phi_\rho(x_k)\le\epsilon$ has to be checked at each inner iteration. However, since the SA algorithm is called to solve subproblems in Algorithm \ref{alg-cons}, we use a more natural termination condition \eqref{epsi-g}. Therefore, $\phi_\rho(x_k)$ is only computed at outer iterations of Algorithm \ref{alg-cons}.
\item[(iii)] Different termination condition for outer iteration is used in Algorithm \ref{alg-cons}. The algorithm in \cite{cgt11} for the deterministic setting is terminated once the criticality measure {$\theta$ at some point} is below some tolerance. However, this cannot be used in Algorithm \ref{alg-cons} for solving {the SNLP problem} \eqref{orig-prob}, because the whole algorithm is a random process, and any specific instance is not sufficient to characterize the performance of criticality measure in average. So we set a maximum iteration number $N$ to terminate the outer iteration of Algorithm \ref{alg-cons}. We will explore the property of the expectation of the output $x_N$ later.
\end{itemize}
\end{rem}

In the following, we shall discuss the $\SFO$-calls complexity of Algorithm \ref{alg-cons}. We assume that the sequence $\{x_k\}$ generated by Algorithm \ref{alg-cons} is bounded. Then {\bf AS.1} indicates that there exist positive constants $\kappa_f$, $\kappa_c$, $\kappa_g$ and $\kappa_J$ such that for all $k$,
\be\label{4-bound}
 f(x_k) \le \kappa_f, \quad  \|c(x_k)\|\le \kappa_c,\quad   \|\nabla f(x_k)\| \le \kappa_g\quad \mbox{and}\quad \|{J}(x_k)\| \le  \kappa_J.
\ee

We first provide an estimate on the optimality of the iterate $x_k$.

\begin{lem}\label{lem3.1}
Let {\bf AS.1} and {\bf AS.4} hold. For fixed $\rho:=\rho_{k-1}$ and any given $\epsilon>0$, if Algorithm \ref{alg-uncons} returns $x_k$ satisfying $\E[\|\tilde{g}_k^r\|^2]\le \epsilon$, then there exists $\lambda_k\in\R^q$ such that
\be \label{exp-res}
\E[\|\nabla f(x_k) + {J}(x_k)^T\lambda_k\|^2] \le 2\epsilon + 2\E[\|G_k - \nabla f(x_k)\|^2],
\ee
where the expectations are taken with respect to the random variables generated in Algorithm \ref{alg-uncons} for solving the ($k$-1)-th subproblem, and $\tilde{g}_k^r$ is defined in \eqref{g-rand_g}.
\end{lem}

{\it Proof.}\quad
Note that the outputs of Algorithm \ref{alg-uncons} are denoted as $x_k = x_{k-1,R_{k-1}}$ and $G_k=G_{k-1,R_{k-1}}$. At the point $x_k$, Algorithm \ref{alg-uncons} generates the next iterate $x_k^+:=x_{k-1,R_{k-1}+1}$ via
\be \label{x_+first}
x_k^+ := \argmin_{u\in\R^n} \left\{ {G_k^T(u-x_k)} + \rho\|c(x_k) + {J}(x_k)(u-x_k)\| + \frac{1}{2\gamma_{k-1,R_{k-1}}}\|u-x_k\|^2\right\}.
\ee
According to the first-order optimality conditions for \eqref{x_+first}, there exists
$
p_k\in \partial \|c(x_k)+{J}(x_k)(x_k^+-x_k)\|
$
such that
\bee
G_k + \rho {J}(x_k)^Tp_k + \frac{1}{\gamma_{k-1,R_{k-1}}}(x_k^+ - x_k) = 0,
\eee
which yields $G_k + \rho {J}(x_k)^Tp_k = \tilde{g}_{k-1,R_{k-1}}^r$. Thus we have the following inequality:
\begin{align}
\|\nabla f(x_k) + \rho {J}(x_k)^Tp_k \|^2 & \le  2\|G_k + \rho {J}(x_k)^Tp_k\|^2 + 2\|G_k - \nabla f(x_k)\|^2 \notag \\
& =  2\|\tilde{g}_{k-1,R_{k-1}}^r\|^2 + 2\|G_k - \nabla f(x_k)\|^2. \label{opt-first}
\end{align}
Hence, by letting $\lambda_k=\rho p_k$ and taking expectation on both sides of \eqref{opt-first}, we obtain \eqref{exp-res}.
\qed

The following lemma shows that, for any given $\epsilon>0$, we can bound $\E[\|\nabla f(x_k) + {J}(x_k)^T\lambda_k\|^2]$ by $\epsilon$ through choosing appropriate total number of $\SFO$-calls and batch sizes when Algorithm \ref{alg-uncons} is applied to solve the NSCO subproblems.
\begin{lem}\label{lem3.2}
Let {\bf AS.1} and {\bf AS.4} hold. For fixed $\rho:=\rho_{k-1}$ and any given $\epsilon>0$, when applying Algorithm \ref{alg-uncons} to minimize $\Phi_\rho$, we choose constant stepsize $\gamma=\gamma_\rho:=1/L_\rho$ and set the total number of $\SFO$-calls $\bar{N}_\rho$ in Algorithm \ref{alg-uncons} as
\be\label{bound-N-cons}
\bar{N}_\rho \ge \max\left\{\frac{\left(4D_{\Phi_\rho}C_2 + 4L_\rho C_3\right)^2}{\epsilon^2} + \frac{128L_\rho D_{\Phi_\rho}}{\epsilon}, \frac{C_1}{L_\rho^2}\right\}.
\ee
where $C_1$, $C_2$ and $C_3$ are defined in \eqref{C_1-C_2-C_3},
\be \label{D-L}
D_{\Phi_\rho} = \Phi_\rho(x_{k-1})-\Phi_\rho^{low}\quad \mbox{and}\quad L_\rho=L_g + \rho L_J.
\ee
We also assume that the batch sizes are chosen to be $m_\rho$:
\be \label{batch-size-m-cons}
m_\rho: = \left\lceil \min \left\{ \bar{N}_\rho, \max\left\{1, \frac{\sigma}{L_\rho}\sqrt{\frac{\bar{N}_\rho}{\tilde{D}}}\right\}\right\}\right\rceil,
\ee
where $\tilde{D}$ is some problem-independent positive constant. Then we have
\be\label{cons-exp-g}
\E[\|\tilde{g}_k^r\|^2] \le \epsilon \quad \mbox{and}\quad \E[\|\tilde{g}_k\|^2]\le \epsilon,
\ee
where the expectations are taken with respect to the random variables generated when the ($k$-1)-th subproblem is solved by Algorithm \ref{alg-uncons}.
Moreover, there exists $\lambda_k\in\R^q$ such that
\be \label{bound-opt}
\E[\|\nabla f(x_k) + {J}(x_k)^T\lambda_k\|^2] \le \epsilon,
\ee
\end{lem}

{\it Proof.}\quad
Let $\epsilon^\prime := \epsilon/4$. Replacing $\epsilon$ by $\epsilon^\prime$ in Theorem \ref{cor2.5}, and using \eqref{bound-N-cons}, we obtain that
\bee
\E[\|\tilde{g}_k^r\|^2] \le \epsilon^\prime \quad \mbox{and}\quad \E[\|\tilde{g}_k\|^2]\le \epsilon^\prime.
\eee
Thus \eqref{cons-exp-g} holds naturally. According to \eqref{delta-bound}, we have $\E[\|G_k - \nabla f(x_k)\|^2] \le \sigma^2/m_\rho$. Similar to Theorem \ref{cor2.5}, we can obtain that
\be \label{exp-G-g}
\E[\|G_k - \nabla f(x_k)\|^2] \le \epsilon^\prime,
\ee
where we have used \eqref{bound-N-cons} and \eqref{batch-size-m-cons}. Therefore, Lemma \ref{lem3.1} indicates
\[\E[\|\nabla f(x_k) + {J}(x_k)^T\lambda_k\|^2] \le 2\epsilon^\prime + 2\epsilon^\prime = \epsilon,\]
i.e., \eqref{bound-opt} holds.
\qed

\begin{rem} \label{rem3.2}
Note that the number of $\SFO$-calls $\bar{N}_\rho$ given in \eqref{bound-N-cons} relies on both $D_{\Phi_\rho}$ and $L_\rho$. Actually both $D_{\Phi_\rho}$ and $L_\rho$ are in the order of $O(\rho)$. To see this, by {\bf AS.1}, we know that for $\rho:=\rho_k$, $k=1,2,\ldots$,
\[
D_{\Phi_\rho} = \Phi_\rho(x_{k-1}) - \Phi_\rho^{low}
               = f(x_{k-1}) + \rho\|c(x_{k-1})\| - \Phi_\rho^{low}                \le \kappa_f + \rho \kappa_c - f^{low},
\]
which implies that $D_{\Phi_\rho}=O(\rho)$. $L_\rho=O(\rho)$ follows directly from \eqref{D-L}.
\end{rem}

Notice that in Algorithm \ref{alg-cons}, for any given $x_k$, $\phi_\rho(x_k)$ plays a key role in adjusting penalty parameters. In the penalty algorithm with exact gradient information proposed by Cartis \etal in \cite{cgt11}, $\phi_{\rho_{k-1}}(x_k)\le\epsilon$ with $G_k$ replaced by $\nabla f(x_k)$ in \eqref{phi} is required as the subproblem termination criterion. However, since an SA algorithm is called to solve subproblems in Algorithm \ref{alg-cons}, a different subproblem termination condition is set to yield \eqref{epsi-g}, namely, $\E[\|\tilde{g}_k^r\|^2]\le\epsilon$. The following lemma provides some interesting relationship between $\E[\|\tilde{g}_k^r\|^2]$ and $\E[\phi_{\rho_{k-1}}(x_k)]$.

\begin{lem}\label{phi-rho}
Let {\bf AS.1} and {\bf AS.4} hold. For fixed $\rho:=\rho_{k-1}\ge1$ and any given $\epsilon>0$, suppose that the iterate $x_k$ is returned by Algorithm \ref{alg-uncons} at the ($k$-1)-th iteration, with stepsizes $\gamma=\gamma_\rho:=1/L_{\rho}$, the number of $\SFO$-calls $\bar{N}_\rho$ satisfying \eqref{bound-N-cons} and batch sizes $m_\rho$ chosen as \eqref{batch-size-m-cons}. Then there exists a positive constant $\bar{C}$ independent of $\rho$ such that
\[
\E[\phi_\rho(x_k)]  \le  2\bar{C}\epsilon^{1/2} + (2\bar{C}L_\rho)^{1/2} \epsilon^{1/4},
\]
where the expectation is taken with respect to random variables generated by Algorithm \ref{alg-uncons} when the ($k$-1)-th subproblem is solved, $\phi_\rho$ is defined in \eqref{phi} and $\bar{C}$ is defined as
\be \label{bar_C}
\bar{C} =\frac{1}{L_J} \kappa_J + \frac{1}{L_g} \left(\kappa_g^2 + 0.25\epsilon \right)^{1/2} ,
\ee
and $L_\rho=L_g+\rho L_J$.
\end{lem}

{\it Proof.}\quad
According to the setting of Algorithm \ref{alg-uncons}, Lemma \ref{lem3.2} shows that $\E[\|\tilde{g}_k^r\|^2]\le \epsilon$.
Recall that starting from $x_k$ Algorithm \ref{alg-uncons} generates the next iterate through
\bee
x_k^+: = \argmin_{u\in\R^n} \left\{ \psi_{\rho,\gamma}(x_k,G_k,u): = {G_k^T(u-x_k)} + \rho \|c(x_k) + {J}(x_k)(u-x_k)\| + \frac{1}{2\gamma}\|u-x_k\|^2 \right\}.
\eee
Then as $\tilde{g}_k^r = (x_k-x_k^+)/\gamma$, we have that
\be \label{E_x_x^+}
\E[\|x_k-x_k^+\|^2] \le \gamma^2\epsilon,
\ee
where the expectation is taken with respect to all the random variables generated by Algorithm \ref{alg-uncons} when the ($k$-1)-th subproblem is solved.

Denote $\Delta\psi_{\rho,\gamma}^k$ as
\[
\Delta\psi_{\rho,\gamma}^k:=\psi_{\rho,\gamma}(x_k,G_k,x_k) - \psi_{\rho,\gamma}(x_k,G_k,x_k^+).
\]
Apparently, $\Delta\psi_{\rho,\gamma}^k>0$. Moreover, it follows from {\bf AS.1} that
\begin{align}
 \Delta\psi_{\rho,\gamma}^k \le &\; \rho\left|\|c(x_k)\| - \|c(x_k)+{J}(x_k)(x_k^+-x_k)\|\right| + \|G_k\|\cdot\|x_k^+-x_k\| - \frac{1}{2\gamma}\|x_k^+-x_k\|^2  \notag \\
\le & \;\rho\kappa_J\|x_k^+-x_k\| + \|G_k\| \cdot\|x_k^+-x_k\| . \label{inq-delta-phi}
\end{align}
For fixed $\rho$, $x_k$ is a random variable generated in the process of Algorithm \ref{alg-uncons}. By taking expectations on both sides of \eqref{inq-delta-phi}, we obtain that
\begin{align*}
\E[\Delta\psi_{\rho,\gamma}^k]&\le  \rho\kappa_J\left(\E[\|x_k^+-x_k\|^2]\right)^{1/2} + \left(\E[\|G_k\|^2]\right)^{1/2}\cdot\left(\E[\|x_k^+-x_k\|^2]\right)^{1/2} \notag  \\
&\le  \rho\gamma\kappa_J\epsilon^{1/2} + \left(\E[\|\nabla f(x_k)\|^2] + \E[\|G_k-\nabla f(x_k)\|^2]\right)^{1/2} \gamma\epsilon^{1/2} 
\\
&\le  \rho\gamma\kappa_J\epsilon^{1/2} + \left(\kappa_g^2 + 0.25\epsilon\right)^{1/2}\gamma\epsilon^{1/2}, \notag
\end{align*}
where the second inequality is from \eqref{E_x_x^+} and the last inequality is due to \eqref{exp-G-g}. According to $\gamma=1/L_\rho$ we have
\begin{align}
\E[\Delta\psi_{\rho,\gamma}^k]
&\le  \left[\frac{1}{L_g+\rho L_J}\rho\kappa_J + \frac{1}{L_g+\rho L_J} \left(\kappa_g^2 + 0.25\epsilon\right)^{1/2}\right]\epsilon^{1/2}\notag \\
& \le  \left[\frac{1}{L_J} \kappa_J + \frac{1}{L_g} \left(\kappa_g^2 + 0.25\epsilon \right)^{1/2} \right]\epsilon^{1/2}  =\bar{C}\epsilon^{1/2}, \label{upp-bound-psi}
\end{align}
where the last inequality is due to $\rho\ge 1$.

We now analyze the property of $\phi_{\rho}(x_k)$, which is defined in \eqref{phi}.
It follows from Lemma 2.5 in \cite{cgt11} that
\bee
\Delta\psi_{\rho,\gamma}^k\ge \frac{1}{2}\min\{1,\gamma\phi_\rho(x_k)\}\phi_\rho(x_k).
\eee
If $1<\gamma\phi_{\rho}(x_k)$, then
\be \label{cond-1}
\phi_\rho(x_k)\le 2\Delta\psi_{\rho,\gamma}^k.
\ee
If $1\ge \gamma\phi_{\rho}(x_k)$, then $\phi_\rho^2(x_k)\le 2\Delta\psi_{\rho,\gamma}^k/\gamma$, which implies
\be \label{cond-2}
\phi_\rho(x_k) \le \gamma^{-1/2}(2\Delta\psi_{\rho,\gamma}^k)^{1/2}.
\ee
Combining \eqref{cond-1} and \eqref{cond-2}, we obtain
\be \label{phi-rho-delta-psi}
\phi_\rho(x_k) \le \max\left\{2\Delta\psi_{\rho,\gamma}^k, \gamma^{-1/2}(2\Delta\psi_{\rho,\gamma}^k)^{1/2}\right\}  \le  2\Delta\psi_{\rho,\gamma}^k + \gamma^{-1/2}(2\Delta\psi_{\rho,\gamma}^k)^{1/2}.
\ee
Taking expectation on both sides of \eqref{phi-rho-delta-psi}, we have
\begin{align*}
\E[\phi_\rho(x_k)] &\le 2\E[\Delta\psi_{\rho,\gamma}^k] + \gamma^{-1/2}\cdot\E[(2\Delta\psi_{\rho,\gamma}^k)^{1/2}]\\
&\le  2\E[\Delta\psi_{\rho,\gamma}^k] +  2^{1/2}\gamma^{-1/2} \cdot(\E[\Delta\psi_{\rho,\gamma}^k])^{1/2}\\
&\le  2\bar{C}\epsilon^{1/2} + (2\bar{C}L_\rho)^{1/2} \epsilon^{1/4},
\end{align*}
where the last inequality is derived from \eqref{upp-bound-psi} and $\gamma=1/L_\rho$.
This completes the proof.
\qed

We next give the main complexity result of Algorithm \ref{alg-cons}.

\begin{thm}\label{thm3.1}
Let {\bf AS.1} and {\bf AS.4} hold. Assume that Algorithm \ref{alg-uncons} is called to solve the NSCO subproblem \eqref{exact-penalty} for fixed $\rho$ at each iteration, with $\gamma=\gamma_\rho:=1/(L_g + \rho L_J)$, the number of $\SFO$-calls $\bar{N}_\rho$ satisfying \eqref{bound-N-cons} and batch sizes $m_\rho$ chosen as \eqref{batch-size-m-cons}. Then Algorithm \ref{alg-cons} returns $x_N$ which satisfies
\be\label{eq1}
\E[\theta(x_N)]  \le\frac{2\bar{C}+(2\bar{C})^{1/2}(L_g+L_J)^{1/2}}{\xi(\rho_0+(N-1)\tau)^{1/2}}\epsilon^{1/4} + \frac{(\kappa_g^2 + 0.25\epsilon)^{1/2}}{(1-\xi)(\rho_0+(N-1)\tau)}
\ee
and
\be\label{eq2}
\E[\|\nabla f(x_N) + J(x_N)^T\lambda_N\|^2]  \le \epsilon, \quad \mbox{for some }\lambda_N\in\R^q,
\ee
where the expectations are taken with respect to all the random variables generated in the process of Algorithm \ref{alg-cons}. Consequently, if we set $N$ as
\be\label{bar-N}
N \ge \hat{N} := \left\lceil\tau^{-1}\tilde{C}\epsilon^{-1/2} -\tau^{-1}\rho_0 + 1\right\rceil,
\ee
where $\tilde{C}=\max\{(4\bar{C}+(8\bar{C})^{1/2}(L_g+L_J)^{1/2})^2\xi^{-2},(4\kappa_g^2 + \epsilon)^{1/2}(1-\xi)^{-1}\}$,
then Algorithm \ref{alg-cons} returns an $\epsilon$-stochastic critical point of \eqref{orig-prob}.

Moreover, Algorithm \ref{alg-cons} finds an $\epsilon$-stochastic critical point of \eqref{orig-prob} after at most
$O(\epsilon^{-3.5})$ $\SFO$-calls.
\end{thm}

{\it Proof.}\quad
Lemma \ref{lem3.2} shows that for any fixed $\rho:=\rho_{k-1}$, $x_k$ returned by Algorithm \ref{alg-uncons} satisfies \eqref{bound-opt}. Because $\rho$ is also a random variable during the process of Algorithm \ref{alg-cons}, \eqref{bound-opt} becomes
\be \label{exp-cond-rho}
\E[\|\nabla f(x_k) + {J}(x_k)^T\lambda_k\|^2|\rho_{[k]}]\le \epsilon,
\ee
where $\rho_{[k]}:=(\rho_1, \ldots, \rho_{k-1})$ and the conditional expectation $\E[\cdot|\rho_{[k]}]$ is taken with respect to the random variables generated by Algorithm \ref{alg-uncons} at the ($k$-1)-th iteration. By further taking expectation with respect to $\rho_{[k]}$ on both sides of \eqref{exp-cond-rho} with $k=N$, we obtain \eqref{eq2}.

We next study the expectation of $\theta(x_N)$, i.e. $\E[\theta(x_N)]$. There are two cases that may happen when Algorithm \ref{alg-cons} terminates, i.e., when $x_N$ is returned as the approximate solution of \eqref{orig-prob}. One case is that $\rho:=\rho_{N-1}$ satisfies \eqref{phi-theta}, namely,
\be \label{upp-bound-theta}
\theta(x_N) \le \frac{\phi_\rho(x_N)}{\xi\rho_{N-1}}.
\ee
The other case is that \eqref{phi-theta} does not hold at $\rho:=\rho_{N-1}$, then it indicates that the inequality $\phi_{\rho} (x_N) < \rho_{N-1}\xi\theta(x_N)$ holds. By \eqref{rho-low-bound} we have
\be\label{the}
\theta(x_N) < \frac{\|G_N\|}{(1-\xi)\rho_{N-1}}.
\ee
Then combining \eqref{upp-bound-theta} and \eqref{the} we obtain
\begin{align}
\theta(x_N) &\le\max\left\{\frac{\phi_\rho(x_N)}{\xi\rho_{N-1}},\frac{\|G_N\|}{(1-\xi)\rho_{N-1}}\right\}\notag\\
&\le \frac{\phi_\rho(x_N)}{\xi\rho_{N-1}} + \frac{\|G_N\|}{(1-\xi)\rho_{N-1}}.\label{est-the}
\end{align}
We first analyze the expectation of $\theta(x_N)$ conditioned on $\rho_{N-1}$, i.e. $\E[\theta(x_N)|\rho_{N-1}]$. In this case, the expectation is taken with respect to the random variables generated when the NSCO subproblem is solved with $\rho=\rho_{N-1}$. On the one hand, Lemma \ref{phi-rho} shows that the expectation of $\phi_{\rho_{N-1}}(x_k)$ satisfies
\bee
\E[(\phi_{\rho_{N-1}}(x_N))|\rho_{N-1}]  \le  2\bar{C}\epsilon^{1/2} + (2\bar{C})^{1/2}(L_g+\rho_{N-1} L_J)^{1/2} \epsilon^{1/4},
\eee
where $\bar{C}$ is defined in \eqref{bar_C}. By taking expectation {on the first term of \eqref{est-the}} conditioned on $\rho_{N-1}$, we have
\begin{align}
\E[\frac{\phi_{\rho_{N-1}}(x_N)}{\xi\rho_{N-1}}|\rho_{N-1}] & \le  \frac{2\bar{C}+(2\bar{C})^{1/2}(L_g+L_J)^{1/2}}{\xi(\rho_{N-1})^{1/2}}\epsilon^{1/4} \label{exp-cond-case1} \\
& \le \frac{2\bar{C}+(2\bar{C})^{1/2}(L_g+L_J)^{1/2}}{\xi(\rho_0+(N-1)\tau)^{1/2}}\epsilon^{1/4}:=E_1, \notag
\end{align}
where the first inequality follows from the facts that $\epsilon \ll 1$ and $\rho_{k}\ge1$ for any $k$, and the second inequality follows from $\rho_{N-1}\ge \rho_0 + (N-1)\tau$. On the other hand, {by taking expectation on the second term of \eqref{est-the}} conditioned on $\rho_{N-1}$, we have
\begin{align}
\E[\frac{\|G_N\|}{(1-\xi)\rho_{N-1}}|\rho_{N-1}] &\le \frac{(\E[\|G_N\|^2|\rho_{N-1}])^{1/2}}{(1-\xi)\rho_{N-1}} \notag \\
&=  \frac{(\E[\|\nabla f(x_k)\|^2|\rho_{N-1}] + \E[\|G_k - \nabla f(x_k)\|^2|\rho_{N-1}])^{1/2}}{(1-\xi)\rho_{N-1}}   \notag \\
&\le  \frac{(\kappa_g^2 + 0.25 \epsilon)^{1/2} }{(1-\xi)\rho_{N-1}} \label{exp-theta-rho1}\\
& \le  \frac{(\kappa_g^2 + 0.25\epsilon)^{1/2}}{(1-\xi)(\rho_0+(N-1)\tau)} :=E_2, \notag
\end{align}
where the second inequality follows from \eqref{exp-G-g} and the last one is due to the fact $\rho_{N-1} \ge \rho_0 + (N-1)\tau$. Then by \eqref{est-the} it yields that $\E[\theta(x_N)|\rho_{N-1}] \le E_1+E_2$. Since both $E_1$ and $E_2$ are fixed constants, after further taking expectation with respect to $\rho_{[N]}$ we obtain $\E[\theta(x_N)]\le E_1+E_2$ which is exactly \eqref{eq1}. Moreover, note that $E_1\le \sqrt{\epsilon}/2$ if
$
N\ge \hat{N}_1:= \lceil \frac{(4\bar{C}+(8\bar{C})^{1/2}(L_g+L_J)^{1/2})^2}{\xi^2\tau}\epsilon^{-1/2} - \frac{\rho_0}{\tau} + 1\rceil
$.
And it gives $E_2 \le \sqrt{\epsilon}/2$ if
$
N\ge \hat{N}_2:=\lceil \frac{(4\kappa_g^2+\epsilon)^{1/2}}{(1-\xi)\tau}\epsilon^{-1/2} - \frac{\rho_0}{\tau} + 1\rceil.
$
Consequently, we have $\E[\theta(x_N)]\le\sqrt{\epsilon}$ if the maximum iteration number $N$ satisfies \eqref{bar-N}.

We now prove the second part of Theorem \ref{thm3.1}. From \eqref{exp-cond-case1} and \eqref{exp-theta-rho1} we know $\E[\theta(x_N)|\rho_{N-1}]\le \sqrt{\epsilon}$, if
\bee
\rho_{N-1} \ge \bar{\rho}:= \tilde{C}\epsilon^{-1/2}.
\eee
Hence, after at most
$
\left\lceil\frac{\bar{\rho} - \rho_0}{\tau}\right\rceil = \hat{N}-1
$
iterations, $\rho_0$ can be increased to no less than $\bar{\rho}$ and we thus have $\E[\theta(x_N)|\rho_{N-1}]\le \sqrt{\epsilon}$.
By taking expectation with respect to $\rho_{[N]}$ we obtain $\E[\theta(x_N)]\le \sqrt{\epsilon}$.
Moreover, from Lemma \ref{lem3.2} we know that for any $k$, to achieve \eqref{exp-cond-rho} at the ($k$-1)-th iteration, Algorithm \ref{alg-uncons} needs at most
$
\max\{\left(4D_{\Phi_\rho}C_2 + 4L_\rho C_3\right)^2\epsilon^{-2} + 128L_\rho D_{\Phi_\rho}\epsilon^{-1}, C_1L^{-2}_\rho\}
$
$\SFO$-calls, where $\rho=\rho_{k-1}$, $D_{\Phi_\rho}=O(\rho)$, $L_\rho=O(\rho)$ and $C_1$, $C_2$, $C_3$ are all constants. Hence, before $\rho$ increases to $\bar{\rho}$, the number of $\SFO$-calls at each iteration is at most in the order of $O(\bar{\rho}^2 \epsilon^{-2})$. Therefore, after at most
\bee
O\left(\hat{N} \bar{\rho}^2 \epsilon^{-2}\right) = O\left(\epsilon^{-3.5}\right)
\eee
$\SFO$-calls, the iterate $x_N$ generated by Algorithm \ref{alg-cons} is an $\epsilon$-stochastic critical point of \eqref{orig-prob}.
\qed


\section{A penalty method with stochastic zeroth-order approximation for SNLP \eqref{orig-prob}}\label{sec:penalty-0th}

In this section, we shall study a penalty method for {SNLP} \eqref{orig-prob}, for which we assume that only noisy function values of $f$ can be obtained via calls to $\SZO$. For any input $x_k$, $\SZO$ outputs a stochastic function value $F(x_k, \xi_k)$, where $\xi_k$ is a random variable whose distribution is supported on $\Xi\subseteq \R^d$ and independent of $x_k$. Furthermore, we assume that $F(x_k,\xi_k)$ is an unbiased estimator of $f(x_k)$. We thus make the following assumption for $\SZO$.
\begin{description}
\item[]\textbf{AS.5}\quad For any $k\ge1$, $F(\cdot,\xi_k)$ is continuously differentiable and $\nabla F(\cdot,\xi_k)$ is Lipschitz continuous with Lipschitz constant $L_g$ for fixed $\xi_k$ and
\begin{align} \label{zero-ord-exp}
\E_{\xi_k}[F(x_k,\xi_k)] = f(x_k).
\end{align}
Throughout this section, we denote
\be \label{G}
G(x_k, \xi_k) =  \nabla_x F(x_k,\xi_k),
\ee
and assume that {\bf AS.4} holds for $G(x_k,\xi_k)$.
\end{description}

As only zeroth-order information of $f$ can be obtained, we need to figure out how to make full use of such information. One of the most popular ways is to apply smoothing techniques. Randomized smoothing techniques have been proposed and fully studied in \cite{dbw12,gl13,glz13,n10}.
We here consider the Gaussian distribution smoothing technique.
For any function $\omega$, given an $n$-dimensional Gaussian random vector $v$, the Gaussian smoothing approximation function of $\omega$ is defined as
\be\label{guassian-smoothing-function}
\omega_\mu(x) := \E_v[\omega(x+\mu v)] =  \frac{1}{(2\pi)^{n/2}}\int \omega(x+\mu v)e^{-\frac{1}{2}\|v\|^2}dv .
\ee
We next cite a lemma which gives some nice properties of the Gaussian smoothing approximate function $\omega_\mu$ in \eqref{guassian-smoothing-function}. This lemma has been proved in \cite{n10} and is also used in \cite{glz13}.

\begin{lem} \label{lem4.1}
If $\omega\in\mathcal{C}_L^{1,1}(\R^n)$ \footnote{$\omega\in\mathcal{C}_L^{1,1}(\R^n)$ means that $\omega: \R^n\to \R$ is continuously differentiable and $\nabla \omega$ is Lipschitz continuous with Lipschitz constant $L$.}, then
\begin{itemize}
\item[a)] $\omega_\mu$ is Lipschitz continuously differentiable with gradient Lipschitz constant $L_\mu\le L$ and
    \[
    \nabla \omega_\mu(x)= \frac{1}{(2\pi)^{n/2}}\int \frac{\omega(x+\mu v) - \omega(x)}{\mu}ve^{-\frac{1}{2}\|v\|^2}dv;
    \]
\item[b)] for any $x\in\R^n$, we have
\begin{align}
& |\omega_\mu(x) - \omega(x)|  \le  \frac{\mu^2}{2}Ln, \label{ine-1}\\
& \|\nabla \omega_\mu(x) - \nabla \omega(x)\|  \le   \frac{\mu}{2}L(n+3)^{\frac{3}{2}}, \label{ine-2} \\
& \E_v\left[\left\|\frac{\omega(x+\mu v) - \omega(x)}{\mu}v\right\|^2\right]  \le
2(n+4)\|\nabla \omega(x)\|^2 + \frac{\mu^2}{2}L^2(n+6)^3; \label{ine-3}
\end{align}
\item[c)] $\omega_\mu$ is convex if $\omega$ is convex.
\end{itemize}
\end{lem}

With the stochastic zeroth-order information of $f$ at $x_k$, namely $F(x_k,\xi_k)$, we can further define the stochastic gradient of $f$ at $x_k$ as
\be \label{G_mu}
G_\mu (x_k,\xi_k,v) = \frac{F(x_k + \mu v, \xi_k) - F(x_k,\xi_k)}{\mu} v.
\ee
From \eqref{zero-ord-exp} and a) of Lemma \ref{lem4.1}, it follows that
\bee
\E_{v,\xi_k}[G_\mu (x_k,\xi_k,v)] =  \nabla f_\mu(x_k).
\eee

When solving \eqref{orig-prob}, the penalty function minimization subproblem in this case is a special NSCO problem in which only noisy function values of $f$ can be obtained via $\SZO$ calls. So we need
to first present {an SA} algorithm, Algorithm \ref{alg-uncons-zero}, with only stochastic zeroth-order information being used for solving NSCO \eqref{nonsmooth-uncons}.

\begin{algorithm}[H]
\caption{\quad \bf Stochastic zeroth-order approximation algorithm for NSCO \eqref{nonsmooth-uncons}}
\label{alg-uncons-zero}
\begin{algorithmic}[1]
\REQUIRE Given $x_1\in\R^n$, maximum iteration number ${N_{in}}$, parameters $\{\gamma_k\}$ with $\gamma_k>0$, batch sizes $\{m_k\}$ with $m_k>0$, a smoothing parameter $\mu>0$. Let $R$ be a random variable following probability distribution $P_R$ which is supported on $\{1,\ldots,{N_{in}}\}$.
\ENSURE {$x_R$.}
\FOR{$k=1,\ldots,R-1$,}
\STATE Call $\SZO$ $m_k$ times to obtain $G_\mu(x_k,\xi_{k,i},v_{k,i})$, $i=1,\ldots,m_k$, where $G_\mu(x_k,\xi_{k,i},v_{k,i})$ is defined in \eqref{G_mu}. Set
    \begin{equation} \label{G_mu_k}
    G_{\mu,k} := \frac{1}{m_k}\sum_{i=1}^{m_k}G_{\mu}(x_k,\xi_{k,i},v_{k,i}).
    \end{equation}
\STATE Compute
\bee
x_{k+1} = \argmin_{u\in\R^n}\psi_{\gamma_k}(x_k,G_{\mu,k},u).
\eee
\ENDFOR
\end{algorithmic}
\end{algorithm}

We denote
\be \label{g-rand_g-zero}
\tilde{g}_{\mu,k} = P_{\gamma_k}(x_k,\nabla f_\mu(x_k)) \quad\mbox{and} \quad \tilde{g}_{\mu,k}^r = P_{\gamma_k}(x_k,G_{\mu,k}),
\ee
where $P_\gamma(x,g)$ is defined in \eqref{proj-gra}. Similar to first-order SA method, we can obtain some properties of Algorithm \ref{alg-uncons-zero}. We next state two main results: Theorems \ref{thm4.1}, \ref{cor4,4} with their proofs given in Appendix A. The following Theorem \ref{thm4.1} provides a bound for $\E[\|\tilde{g}_{\mu,R}^r\|^2]$.

\begin{thm} \label{thm4.1}
Let {\bf AS.1-5} hold. Suppose that the stepsizes $\{\gamma_k\}$ in Algorithm \ref{alg-uncons-zero} are chosen such that $0<\gamma_k\le 2/L$, $k=1,\ldots,N$, with $\gamma_k<2/L$ for at least one $k$, where $L=L_g+L_hL_J$. Moreover, suppose that the probability mass function $P_R$ is chosen as in \eqref{prob_R}, and suppose that there exists $\kappa_g>0$ such that $\|\nabla f(x_k)\|\le\kappa_g$ for any $k$.
Then for any $N\ge1$, we have
\be\label{proj-gra-rand-bound-zero}
\E[\|\tilde{g}_{\mu,R}^r\|^2]\le \frac{D_{\Phi_h} + \mu^2 L_g n + \tilde{\sigma}^2\sum_{k=1}^{N_{in}}(\gamma_k/m_k)}{\sum_{k=1}^{N_{in}}(\gamma_k - L\gamma_k^2/2)},
\ee
where the expectation is taken with respect to $R$, $\xi_{[{N_{in}}]}:=(\xi_1,\ldots,\xi_{N_{in}})$ with $\xi_k:=(\xi_{k,1},\ldots,\xi_{k,m_k})$ and $v_{[{N_{in}}]}:=(v_1,\ldots,v_{N_{in}})$ with $v_k:=(v_{k,1},\ldots,v_{k,m_k})$, and $D_{\Phi_h}$ is defined in \eqref{D_Phih} and $\tilde{\sigma}^2$ is defined as
\bee
\tilde{\sigma}^2 =  2(n+4)[\kappa_g^2 + \sigma^2 + \mu^2L_g^2(n+4)^2].
\eee
\end{thm}

By specializing the settings of Algorithm \ref{alg-uncons-zero}, we obtain the following complexity result.
\begin{thm}\label{cor4,4}
Let assumptions {\bf AS.1-5} hold. Suppose that in Algorithm \ref{alg-uncons-zero}, $\gamma_k=1/L$ where $L=L_g + L_h L_J$, the probability mass function $P_R$ is chosen as \eqref{prob_R}, and there exists $\kappa_g>0$ such that $\|\nabla f(x_k)\|\le\kappa_g$ for all $k$. {Denote $\bar{N}$ as the total number of $\mathcal{SZO}$-calls in Algorithm \ref{alg-uncons-zero}. For any given constant $\epsilon>0$, suppose that $\bar{N}$ satisfies}
\be \label{bar_N-zero}
\bar{N} \ge  \max\left\{\frac{(16D_{\Phi_h}/\sqrt{\tilde{D}_2} + L\tilde{C}_1)^2}{\epsilon^2} + \frac{{112}LL_g\tilde{D}_1(n+4)+64LD_{\Phi_h}}{\epsilon},\frac{1}{L^2\tilde{D}_2}\right\},
\ee
where $\tilde{D}_1, \tilde{D}_2$ are two problem-independent positive constants and
\be \label{tilde_C}
\tilde{C}_1 = 24(n+4)(\kappa_g^2+\sigma^2)\sqrt{\tilde{D}}_2.
\ee
Suppose that the smoothing parameter $\mu$ satisfies
\be \label{mu}
\mu \le  \sqrt{\frac{\tilde{D}_1}{\bar{N}}},
\ee
and the batch sizes $m_k=m$ satisfy
\be \label{m}
m =  \left\lceil \min\left\{ \bar{N} , \max\left\{1, \frac{1}{L}\cdot\sqrt{\frac{\bar{N}}{\tilde{D}_2}} \right\} \right\} \right\rceil,
\ee
Then we have
\be\label{eps-exp}
\E[\|\tilde{g}_{\mu,R}^r\|^2]\le\epsilon \quad \mbox{and}\quad \E[\|\tilde{g}_R\|^2]\le\epsilon,
\ee
where the expectations are taken with respect to $R$, $\xi_{[{N_{in}}]}$ and $v_{[{N_{in}}]}$. $\tilde{g}_{k}$ and $\tilde{g}_{\mu,k}^r$ are defined in \eqref{g-rand_g} and \eqref{g-rand_g-zero} respectively. Thus, it follows that the number of $\SZO$-calls required by Algorithm \ref{alg-uncons-zero} to achieve $\E[\|\tilde{g}_{\mu,R}^r\|^2]\le\epsilon$ and $\E[\|\tilde{g}_R\|^2]\le\epsilon$ is in the order of $O(\epsilon^{-2})$.
\end{thm}



We are now ready to present a stochastic zeroth-order penalty method for solving \eqref{orig-prob}. In each iteration, Algorithm \ref{alg-uncons-zero} is called to minimize the penalty function. The strategy to update penalty parameters is the same as the one applied in Algorithm \ref{alg-cons}.

\begin{algorithm}[H]
\caption{\quad \bf Penalty method with stochastic zeroth-order approximation for \eqref{orig-prob}}
\label{alg-pen-zero}
\begin{algorithmic}[1]
\REQUIRE Given maximum iteration number $N$, tolerance $\epsilon\in(0,1)$, initial smoothing parameter $\mu_0$, steering parameter $\xi\in(0,1)$, initial iterate $x_1\in\R^n$, $G_{\mu_0}^1\in\R^n$, penalty parameter $\rho_0\ge1$ and minimal increase factor $\tau>0$. Set $k:=1$.
\ENSURE {$x_N$.}
\FOR{$k=1,\ldots,N-1$}
\STATE {Step (a):} Find $\rho:=\rho_k\ge\rho_{k-1}+\tau$ satisfying
\bee 
\phi_{\rho,\mu_{k-1}}(x_k) \ge \rho \xi\theta(x_k),
\eee
where $\theta(x)$ is defined in \eqref{theta} and
\be\label{phi-mu}
\phi_{\rho,\mu_{k-1}}(x_k) = \rho \|c(x_k)\| - \min_{\|s\|\le1} \{ \langle G_{\mu_{k-1}}^k, s \rangle + \rho\|c(x_k) + {J}(x_k)s\| \},
\ee
\STATE {Step (b)}: Apply Algorithm \ref{alg-uncons-zero} with smoothing parameter $\mu_k$, initial iterate $x_{\mu_k,1}:=x_{k}$ {and $\bar{N}_\rho$ $\SZO$-calls} to solve the subproblem
\bee
\min_{x\in\R^n}\quad \Phi_{\rho_k}(x) = f(x) + \rho_k\|c(x)\|.
\eee
returning $x_{k+1}:=x_{\mu_k,R_k}$ and $G_{\mu_k}^{k+1}:=G_{\mu_k,R_k}$, for which
\bee 
\E[\|\tilde{g}_{\mu_k,R_k}^r\|^2] \le \epsilon,
\eee
where ``$x_{\mu_k,R_k}$" denotes the $R_k$-th iterate generated by Algorithm \ref{alg-uncons-zero} with smoothing parameter $\mu_k$ when solving the $k$-th subproblem and $\tilde{g}_{\mu,k}^r$ is defined in \eqref{g-rand_g-zero}, and the expectation is taken with respect to the random variables generated in this inner iteration.
\ENDFOR
\end{algorithmic}
\end{algorithm}

{Similar to the arguments in Remark \ref{rem3.1}, Step (a) in Algorithm \ref{alg-pen-zero} is well-defined.} Assume that the sequence of iterates $\{x_k\}$ generated by Algorithm \ref{alg-pen-zero} is bounded. Then {\bf AS.1} indicates that {there exist
positive constants $\kappa_f, \kappa_c, \kappa_g$ and $\kappa_J$ such that }  
\eqref{4-bound} holds {for all $k$}.

In the following lemma, we provide a measure on the optimality of each iterate $x_k$.

\begin{lem}
Let assumptions {\bf AS.1} and {\bf AS.4-5} hold. For fixed $\rho:=\rho_{k-1}$ and any given positive constant $\epsilon$, if $x_k$ satisfies that $\E[\|\tilde{g}_{\mu_{k-1},R_{k-1}}^r\|^2]\le \epsilon$, then there exists $\lambda_k\in\R^q$ such that
\be \label{exp-res-zero}
\E[\|\nabla f(x_k) + {J}(x_k)^T\lambda_k\|^2]  \le 4\|G_{\mu_{k-1}}^k - \nabla f_{\mu_{k-1}}(x_k)\|^2 + \mu_{k-1}^2L_g^2(n+3)^3 + 2\epsilon,
\ee
where the expectation is taken with respect to the random variables generated by Algorithm \ref{alg-uncons-zero} when the ($k$-1)-th subproblem is solved, and $\tilde{g}_{\mu,k}^r$ is defined in \eqref{g-rand_g-zero}.
\end{lem}

{\it Proof.}\quad
By the construction of Algorithm \ref{alg-pen-zero}, $G_{\mu_{k-1}}^{k} = G_{\mu_{k-1},R_{k-1}}$ for some $R_{k-1}$. At the iterate $x_k$, Algorithm \ref{alg-uncons-zero} generates the next point $x_k^+$ through
\be \label{x_k-x_+}
x_k^+ := \mbox{arg}\min_{u\in\R^n} \left\{ {(G_{\mu_{k-1}}^k)^T(u- x_k)} + \rho\left\|c(x_k) + {J}(x_k)(u-x_k)\right\| + \frac{1}{2\gamma_{k-1,R_{k-1}}}\|u-x_k\|^2\right\}.
\ee
From the first-order optimality conditions for \eqref{x_k-x_+}, we know that there exists
$p_k\in \partial\|c(x_k)+{J}(x_k)(x_k^+-x_k)\|$
such that
\[
G_{\mu_{k-1}}^k + \rho {J}(x_k)^T p + \frac{1}{\gamma_{k-1,R_{k-1}}}(x_k^+-x_k) = 0,
\]
which shows $G_{\mu_{k-1}}^k + \rho {J}(x_k)^T p_k = \tilde{g}_{\mu_{k-1},R_{k-1}}^r$. Hence we have
\begin{align}
 \|\nabla f(x_k) + \rho {J}(x_k)^T p_k\|^2
 \le & \; 2\|G_{\mu_{k-1}}^k - \nabla f(x_k)\|^2 + 2\|G_{\mu_{k-1}}^k + \rho {J}(x_k)^T p_k\|^2  \notag  \\
 = & \;2\|G_{\mu_{k-1}}^k - \nabla f(x_k)\|^2 + 2 \|\tilde{g}_{\mu_{k-1},R_{k-1}}^r\|^2 \notag  \\
 \le & \;4\|G_{\mu_{k-1}}^k - \nabla f_{\mu_{k-1}}(x_k)\|^2 + 4 \|\nabla f_{\mu_{k-1}}(x_k) - \nabla f(x_k)\|^2 + 2 \|\tilde{g}_{\mu_{k-1},R_{k-1}}^r\|^2  \notag \\
 \le & \;4\|G_{\mu_{k-1}}^k - \nabla f_{\mu_{k-1}}(x_k)\|^2 + \mu_{k-1}^2L_g^2(n+3)^3 + 2 \|\tilde{g}_{\mu_{k-1},R_{k-1}}^r\|^2, \label{opt-zero}
\end{align}
where the last inequality follows from \eqref{ine-2}. Therefore, by taking expectation on both sides of \eqref{opt-zero} with respect to the random variables generated by Algorithm \ref{alg-uncons-zero} when solving the ($k$-1)-th subproblem, we obtain \eqref{exp-res-zero} by letting $\lambda_k=\rho p_k$.
\qed

We show in the following lemma that for any given positive constant $\epsilon$, we can bound $\E[\|\nabla f(x_k) + {J}(x_k)^T\lambda_k\|^2]$ by $\epsilon$ through choosing appropriate total number of $\SZO$ calls $\bar{N}$, the batch size $m$ and the smoothing parameter $\mu$ at each iteration for any fixed $\rho=\rho_{k-1}$.

\begin{lem}\label{lem5.2}
Let {\bf AS.1} and {\bf AS.4-5} hold. For fixed $\rho:=\rho_{k-1}$ and any given positive constant $\epsilon$, suppose that when applying Algorithm \ref{alg-uncons-zero} to minimize $\Phi_\rho$, we choose the constant stepsizes $\gamma_k=\gamma_\rho:=1/L_\rho$ and the total number of $\SZO$-calls $\bar{N}_\rho$ satisfies
\be\label{bound-N-cons-zero}
\bar{N}_\rho \ge  \max\left\{\frac{(64D_{\Phi_\rho}/\sqrt{\tilde{D}_2} + 4L_\rho\tilde{C}_1)^2}{\epsilon^2} + \frac{448L_\rho L_g\tilde{D}_1(n+4) + 256L_\rho D_{\Phi_\rho}}{\epsilon},\frac{1}{L_\rho^2\tilde{D}_2}\right\},
\ee
where $D_{\Phi_\rho}$ and $L_\rho$ are defined in \eqref{D-L},
$\tilde{C}_1$ is defined in \eqref{tilde_C}, and $\tilde{D}_1$ and $\tilde{D}_2$ are two problem-independent positive scalars. Also suppose that the batch sizes are chosen equal to $m_\rho$ defined as
\be \label{batch-size-m-cons-zero}
m_\rho : =  \left\lceil \min\left\{ \bar{N}_\rho , \max\left\{1, \frac{1}{L_\rho}\cdot\sqrt{\frac{\bar{N}_\rho}{\tilde{D}_2}} \right\} \right\} \right\rceil.
\ee
Besides, the smoothing parameter $\mu_{k-1}$ is assumed to satisfy
\begin{equation}\label{mu_k}
\mu_{k-1}\le\sqrt{\frac{\tilde{D}_1}{\bar{N}_\rho}}.
\end{equation}
Then for $x_k:=x_{k-1,R_{k-1}}$ we have
\be\label{cons-exp-g-zero}
\E[\|\tilde{g}_{\mu_{k-1},R_{k-1}}^r\|^2] \le \epsilon, \qquad \E[\|\tilde{g}_k\|^2]\le \epsilon,
\ee
and there exists $\lambda_k\in\R^q$ such that
\be \label{bound-opt-zero}
\E[\|\nabla f(x_k) + {J}(x_k)^T\lambda_k\|^2] \le \epsilon,
\ee
where the expectations are taken with respect to all the random variables generated when the ($k$-1)-th subproblem being solved. 
\end{lem}

{\it Proof.}\quad
First, by letting $\epsilon^\prime = \epsilon/4$, similar to the analysis in Theorem \ref{cor4,4} by replacing $\epsilon$ with $\epsilon^\prime$, we can prove that the choice of $\bar{N}_\rho$ in \eqref{bound-N-cons-zero} can ensure that
$
\E[\|\tilde{g}_{\mu_{k-1},R_{k-1}}^r\|^2] \le \epsilon^\prime$ and $\E[\|\tilde{g}_k\|^2]\le \epsilon^\prime.
$
Therefore, \eqref{cons-exp-g-zero} holds naturally.

Second, noticing that $x_k=x_{\mu_{k-1},R_{k-1}}$ and $G_{\mu_{k-1}}^k = G_{\mu_{k-1},R_{k-1}}$, by \eqref{G-mu-f-zero} we have
\be \label{diff-mu}
\E[\|G_{\mu_{k-1}}^k - \nabla f_{\mu_{k-1}}(x_k)\|^2]  \le  \frac{\tilde{\sigma}_{k-1}^2}{m_\rho},
\ee
where the expectation is taken with respect to the random variables generated by Algorithm \ref{alg-uncons-zero}, and
\be \label{tilde-sigma-k}
\tilde{\sigma}_{k-1} = 2(n+4)[\kappa_g^2 + \sigma^2 + \mu^2_{k-1}L_g^2(n+4)^2].
\ee
So \eqref{exp-res-zero} implies that
\be  \label{exp-opt-1}
\E[\|\nabla f(x_k) + {J}(x_k)^T\lambda_k\|^2]  \le \frac{4\tilde{\sigma}_{k-1}^2}{m_\rho} + \mu_{k-1}^2L_g^2(n+3)^3 + 2\epsilon^\prime.
\ee
Let us consider the first two terms on the right hand side of \eqref{exp-opt-1}. According to the definition of $\tilde{\sigma}_{k-1}$ in \eqref{tilde-sigma-k} and the choice of $\mu_{k-1}$ satisfying \eqref{mu_k}, we have
\begin{align*}
\frac{4\tilde{\sigma}_{k-1}^2}{m_\rho} + \mu_{k-1}^2L_g^2(n+3)^3
 \le\; & \frac{4}{3}\left[ \frac{3\tilde{\sigma}_{k-1}^2}{m_\rho} + \mu_{k-1}^2L_g^2(n+3)^3\right]\\
 \le\; & \frac{4}{3}\left[\frac{6(n+4)(k_g^2+\sigma^2)}{m_\rho} + \frac{6(n+4)^3L_g^2}{m_\rho}\cdot\frac{\tilde{D}_1}{\bar{N}} + \frac{\tilde{D}_1}{\bar{N}}\cdot L_g^2(n+3)^3\right] \\
 \le\; & \frac{4}{3}\zeta :=\frac{4}{3}\left[\frac{6(n+4)(k_g^2+\sigma^2)}{m_\rho} + \frac{7L_g^2\tilde{D}_1(n+4)^3}{\bar{N}}\right].
\end{align*}
Note that $\zeta$ is less than the right hand side of \eqref{right-exp-R}. Following the analysis in Theorems \ref{cor2.5} and \ref{cor4,4} we obtain that the choice of $\bar{N}_\rho$ and $m_\rho$ in \eqref{bound-N-cons-zero} and \eqref{batch-size-m-cons-zero} can ensure
\be \label{comp}
\frac{4\tilde{\sigma}_{k-1}^2}{m} + \mu_{k-1}^2L_g^2(n+3)^3  \le  \frac{4}{3}\cdot\epsilon^\prime < 2\epsilon^\prime.
\ee
Combining \eqref{exp-opt-1} and \eqref{comp} gives \eqref{bound-opt-zero}.
\qed

\begin{rem}
Note that in Lemma \ref{lem5.2}, the number of $\SZO$-calls $\bar{N}_\rho$ in \eqref{bound-N-cons-zero} depends on both $L_\rho$ and $D_{\Phi_\rho}$. Similar to the analysis in Remark \ref{rem3.2}, we obtain that $D_{\Phi_\rho}=O(\rho)$ and $L_\rho=O(\rho)$. Since $\tilde{C}_1$, $\tilde{D}_1$ and $\tilde{D}_2$ are all constants independent with $\rho$, $\bar{N}_\rho$ is in the order of $O(\rho^2\epsilon^{-2})$.
\end{rem}

Analogous to Lemma \ref{phi-rho}, we give an estimate of $\E[\phi_{\rho,\mu_{k-1}}(x_k)]$ in the following lemma.

\begin{lem}
Let {\bf AS.1} and {\bf AS.4-5} hold. For fixed $\rho=\rho_{k-1}$ and any given positive constant $\epsilon$, suppose that the iterate $x_k$ is returned by Algorithm \ref{alg-uncons-zero} at the ($k$-1)-th iteration with the same settings as in Lemma \ref{lem5.2}. Then we have
\bee 
\E[\phi_{\rho,\mu_{k-1}}(x_k)]  \le  2\bar{C}\epsilon^{1/2} + (2\bar{C})^{1/2}(L_g+\rho L_J)^{1/2}\epsilon^{1/4},
\eee
where the expectation is taken with respect to random variables generated by Algorithm \ref{alg-uncons-zero} when solving the ($k$-1)-th subproblem, $\phi_{\rho,\mu_{k-1}}$ is defined in \eqref{phi-mu} and $\bar{C}$ is defined in \eqref{bar_C}.
\end{lem}

{\it Proof.}\quad
The idea of the proof is similar to Lemma \ref{phi-rho}. We only need to estimate $\E[\|G_{\mu_{k-1}}^k - \nabla f(x_k)\|^2]$. By \eqref{diff-mu} and \eqref{ine-2}, we have
\begin{align*}
\E[\|G_{\mu_{k-1}}^k - \nabla f(x_k)\|^2]& \le 2\E[\|G_{\mu_{k-1}}^k - \nabla f_{\mu_{k-1}}(x_k)\|^2] + 2 \E[\|\nabla f_{\mu_{k-1}}(x_k) - \nabla f(x_k)\|^2]  \\
& \le  \frac{2\tilde{\sigma}_{k-1}^2}{m_\rho} + \frac{1}{2}\mu_{k-1}^2L_g^2(n+3)^3  < \epsilon^\prime = \frac{1}{4}\epsilon,
\end{align*}
where the last inequality follows from \eqref{comp}. The rest of the proof is the same as Lemma \ref{phi-rho}.
\qed

We now conclude this section by giving the main result on the total $\SZO$-calls worst-case complexity for Algorithm \ref{alg-pen-zero}. The proof is essentially the same as Theorem \ref{thm3.1}, so we only state the result and omit the proof.

\begin{thm}
Let {\bf AS.1} and {\bf AS.4-5} hold. Assume that Algorithm \ref{alg-uncons-zero} is applied to solve the stochastic subproblem \eqref{exact-penalty} for fixed $\rho$ at each iteration, with $\gamma=\gamma_\rho:=1/(L_g+\rho L_J)$, the number of $\SZO$-calls $\bar{N}_\rho$ satisfying \eqref{bound-N-cons-zero}, batch sizes $m_\rho$ chosen as \eqref{batch-size-m-cons-zero}, and smoothing parameters satisfying \eqref{mu_k}. Then Algorithm \ref{alg-pen-zero} either returns an $\epsilon$-stochastic critical point of \eqref{orig-prob}, or returns $x_N$ which satisfies
\[
\E[\theta(x_N)]  \le\frac{2\bar{C}+(2\bar{C})^{1/2}(L_g+L_J)^{1/2}}{\xi(\rho_0+(N-1)\tau)^{1/2}}\epsilon^{1/4} + \frac{(\kappa_g^2 + 0.25\epsilon)^{1/2}}{(1-\xi)(\rho_0+(N-1)\tau)}
\]
and
\[
\E[\|\nabla f(x_N) + J(x_N)^T\lambda_N\|^2]  \le \epsilon, \quad \mbox{for some }\lambda_N\in\R^q,
\]
where the expectations are taken with respect to all the random variables generated in the process of Algorithm \ref{alg-pen-zero}. Consequently, if we set $N$ satisfying \eqref{bar-N},
then Algorithm \ref{alg-pen-zero} must return an $\epsilon$-stochastic critical point of \eqref{orig-prob}. Moreover, Algorithm \ref{alg-pen-zero} can always find an $\epsilon$-stochastic critical point of \eqref{orig-prob} after at most $O(\epsilon^{-3.5})$ $\SZO$-calls.
\end{thm}

\section{Conclusions}\label{sec:conclusions}

In this paper, we proposed a class of penalty methods with stochastic approximation for solving stochastic nonlinear programming problems. We assumed that only the first-order or zeroth-order information of the objective function was available via subsequent calls to a stochastic first-order or zeroth-order oracle. In each iteration of the penalty methods, we minimized a nonconvex and nonsmooth penalty function to update the iterate. The worst-case complexity of calls to the stochastic first-order (or zeroth-order) oracle for the proposed penalty methods for obtaining an $\epsilon$-stochastic critical point was analyzed.

\section*{Acknowledgements}
We would like to thank two anonymous referees for their insightful comments and suggestions that have helped us improve the presentation of this paper greatly.


\bibliographystyle{plain}
\bibliography{All}

\appendix \sectionfont{Appendix }

\section{Proofs of Theorems \ref{thm2.2}, \ref{cor2.5}, \ref{thm4.1} and \ref{cor4,4}}\label{appA}

In this appendix, we give the detailed proofs of Theorems \ref{thm2.2}, \ref{cor2.5}, \ref{thm4.1} and \ref{cor4,4}. First, we need to prepare some lemmas.

%
%

The following lemma provides a bound for the size of $P_\gamma (x,g)$ defined in \eqref{proj-gra}.

\begin{lem}\label{lem2.2}
Let \textbf{AS.1-2} hold and $P_\gamma (x,g)$ be defined in \eqref{proj-gra}. Then for any $x\in \R^n$, $g\in\R^n$ and $\gamma>0$, we have
\begin{equation}\label{g-proj-gra}
{g^T P_\gamma (x,g) } \ge  \left(1-\frac{1}{2}\gamma L_hL_J \right)\left\|P_\gamma(x,g)\right\|^2 + \frac{1}{\gamma}\left[h(c(x^+))-h(c(x))\right].
\end{equation}
\end{lem}

{\it Proof.}
From the optimality conditions for \eqref{x_+}, it follows that there exists
$
p\in\partial h(c(x)+{J}(x)(x^+-x))
$ such that
$
{( g+{J}(x)^Tp + \frac{1}{\gamma}(x^+-x) )^T( u-x^+) }\ge 0,$ for any $u\in \R^n.
$
Specifically, by letting $u=x$ we obtain
\begin{align*}
{g^T (x - x^+)} &\ge \frac{1}{\gamma}\|x^+-x\|^2 + { p^TJ(x)(x^+-x)} \rangle \ge \frac{1}{\gamma}\|x^+-x\|^2 + h(c(x)+{J}(x)(x^+-x))-h(c(x)),
\end{align*}
where the second inequality is due to the convexity of $h$. {\bf AS.1-2} implies that
\begin{align*}
 |h(c(x^+)) - h(c(x)+{J}(x)(x^+-x))|
 \le\; & L_h\|c(x^+)-(c(x)+{J}(x)(x^+-x))\|  \\
 \le\; & L_h\left\|\int_0^1 [J(x+t(x^+-x))-{J}(x)](x^+-x)dt\right\| \\
 =\; & \frac{1}{2}L_hL_J\|x^+-x\|^2.
\end{align*}
We thus obtain the following bound for $\langle g, x-x^+ \rangle$:
\bee
{g^T(x-x^+ )} \ge \left(\frac{1}{\gamma}-\frac{1}{2}L_hL_J\right)\|x^+-x\|^2 + h(c(x^+))-h(c(x)).
\eee
Therefore, \eqref{g-proj-gra} follows from the definition of $P_\gamma(x,g)$ in \eqref{proj-gra}.
\qed

The following lemma shows that $P_\gamma(x,g)$ is Lipschitz continuous with respect to $g$.
\begin{lem}\label{cor2.1}
Let {\bf AS.1-2} hold and $P_\gamma(x,g)$ be defined in \eqref{proj-gra}. Then for any $g_1,g_2\in\R^n$, we have
\bee
\|P_\gamma(x,g_1)-P_\gamma(x,g_2)\| \le \|g_1 - g_2\|.
\eee
\end{lem}

{\it Proof.} According to \eqref{proj-gra}, letting $x_1^+$ and $x_2^+$ be given through \eqref{x_+} with $g$ replaced by $g_1$ and $g_2$, it suffices to prove that
 $
 \|x_1^+-x_2^+\|  \le \gamma\|g_1-g_2\|.
 $
From the optimality conditions for \eqref{x_+}, there exist
$p_1\in\partial h(c(x)+{J}(x)(x_1^+-x))$ and $p_2\in\partial h(c(x)+{J}(x)(x_2^+-x))$
such that the following two equalities hold:
\begin{align}
( g_1 + {J}(x)^Tp_1 + \frac{1}{\gamma}(x_1^+ - x))^T( u-x_1^+ )&\ge 0, \quad \forall u\in\R^n, \label{opt-1}\\
( g_2 + {J}(x)^Tp_2 + \frac{1}{\gamma}(x_2^+ - x))^T( u-x_2^+ )&\ge 0, \quad \forall u\in\R^n. \label{opt-2}
\end{align}
Letting $u=x_2^+$ in \eqref{opt-1} and using the fact that $h$ is convex, we have
\begin{align}
{ g_1^T (x_2^+-x_1^+)}\ge &\;\frac{1}{\gamma} (x-x_1^+)^T( x_2^+-x_1^+) + p_1^T {J}(x)(x_1^+ - x_2^+) \notag\\
\ge &\; \frac{1}{\gamma} (x-x_1^+)^T( x_2^+-x_1^+) + h(c(x)+{J}(x)(x_1^+-x)) - h(c(x)+{J}(x)(x_2^+-x)). \label{ineq-1}
\end{align}
Similarly, letting $u=x_1^+$ in \eqref{opt-2} we obtain
\be\label{ineq-2}
g_2^T( x_1^+-x_2^+)
\ge \frac{1}{\gamma} (x-x_2^+)^T( x_1^+-x_2^+) + h(c(x)+{J}(x)(x_2^+-x)) - h(c(x)+{J}(x)(x_1^+-x)).
\ee
Summing up \eqref{ineq-1} and \eqref{ineq-2}, we obtain
\begin{equation*}
\|g_1-g_2\|\|x_1^+-x_2^+\| \ge ( g_1-g_2)^T( x_2^+-x_1^+) \ge \frac{1}{\gamma}\|x_1^+ - x_2^+\|^2,
\end{equation*}
which completes the proof.
\qed

We now give the proof of Theorem \ref{thm2.2}.

{\it Proof of Theorem \ref{thm2.2}.}
Denote $\delta_k := G_k - \nabla f(x_k)$. From {\bf AS.1}, we have
\begin{align*}
f(x_{k+1}) & \le  f(x_k) + \nabla f(x_k)^T( x_{k+1} - x_k) + \frac{L_g}{2} \|x_{k+1} - x_k\|^2\\
           & =  f(x_k) + G_k^T( x_{k+1}-x_k) + \frac{L_g}{2} \|x_{k+1} - x_k\|^2 - \langle \delta_k, x_{k+1} - x_k \rangle.
\end{align*}
From the definition of $x_{k+1}$ in \eqref{x_k+1}, it follows that $x_k - x_{k+1} = \gamma_k \tilde{g}_k^r$. According to Lemma \ref{lem2.2} with $g$ replaced by $G_k$ and $x=x_k$ and $\gamma=\gamma_k$, we obtain
\bee
f(x_{k+1}) \le f(x_k) - \left(\gamma_k - \frac{L}{2}\gamma_k^2\right)\|\tilde{g}_k^r\|^2 - h(c(x_{k+1})) + h(c(x_k)) + \gamma_k \delta_k^T\tilde{g}_k^r,
\eee
which implies that
\bee
\Phi_h(x_{k+1}) \le \Phi_h(x_k) - \left( \gamma_k - \frac{L}{2}\gamma_k^2\right)\|\tilde{g}_k^r\|^2 + \gamma_k \delta_k^T\tilde{g}_k + \gamma_k \delta_k^T( \tilde{g}_k^r - \tilde{g}_k).
\eee
Note that it follows from Lemma \ref{cor2.1} with $g_1 = G_k$ and $g_2 = \nabla f(x_k)$ that
\bee
\delta_k^T(\tilde{g}_k^r - \tilde{g}_k ) \le \|\delta_k\|\|\tilde{g}_k^r - \tilde{g}_k\|\le  \|\delta_k\|\|G_k - \nabla f(x_k)\| = \|\delta_k\|^2.
\eee
It yields that
\be \label{Phi-k+1-k}
\Phi_h(x_{k+1}) \le \Phi_h(x_k) - \left( \gamma_k - \frac{L}{2}\gamma_k^2\right)\|\tilde{g}_k^r\|^2 + \gamma_k \delta_k^T\tilde{g}_k +  \gamma_k\|\delta_k\|^2.
\ee
Summing up \eqref{Phi-k+1-k} for $k=1,\ldots,N_{in}$ and noticing that $\gamma_k\le 2/L$, we have
\begin{align}
\sum_{k=1}^{N_{in}} \left( \gamma_k - \frac{L}{2}\gamma_k^2\right) \|\tilde{g}_k^r\|^2 &\le \Phi_h(x_1) - \Phi_h(x_{N_{in}+1}) + \sum_{k=1}^{N_{in}} \{\gamma_k\delta_k^T\tilde{g}_k +  \gamma_k \|\delta_k\|^2\} \notag \\
& \le  \Phi_h(x_1) - \Phi_h^{low} + \sum_{k=1}^{N_{in}} \{\gamma_k\delta_k^T\tilde{g}_k +  \gamma_k \|\delta_k\|^2\}. \label{sum-g-r}
\end{align}
Notice that $x_k$ is a random variable as it is a function of $\xi_{[k-1]}$, generated in the algorithm process. By {\bf AS.4} we have $\E[\delta_k^T\tilde{g}_k|\xi_{[k-1]}]=0$ and
\be \label{delta-bound}
\E[\|G_k - \nabla f(x_k)\|^2]=\E[\|\delta_k\|^2] = \frac{1}{m_k^2}\sum_{i=1}^{m_k}\E[\|\delta_{k,i}\|^2]\le\frac{\sigma^2}{m_k},
\ee
where $\delta_{k,i}=G(x_k, \xi_{k,i}) - \nabla f(x_k)$. Taking the expectation on both sides of \eqref{sum-g-r} with respect to $\xi_{[{N_{in}}]}$, we obtain that
\bee
\sum_{k=1}^{N_{in}} \left( \gamma_k - \frac{L}{2}\gamma_k^2 \right)\E_{\xi_{[{N_{in}}]}}[\|\tilde{g}_k^r\|^2] \le \Phi_h(x_1) - \Phi_h^{low} + \sigma^2 \sum_{k=1}^{N_{in}}\frac{\gamma_k}{m_k}.
\eee
Since $R$ is a random variable with probability mass function $P_R$, it follows that
\bee
\E[\|\tilde{g}_R^r\|^2] = \E_{R,\xi_{[{N_{in}}]}}[\|\tilde{g}_R^r\|^2]=\frac{\sum_{k=1}^{N_{in}}\left( \gamma_k - L\gamma_k^2/2\right)\E_{\xi_{[{N_{in}}]}}[\|\tilde{g}_k^r\|^2]}{\sum_{k=1}^{N_{in}}\left( \gamma_k - L\gamma_k^2/2\right)},
\eee
which proves \eqref{proj-gra-rand-bound}.
\qed

Following from Theorem \ref{thm2.2}, we now prove Theorem \ref{cor2.5}.

{\it Proof of Theorem \ref{cor2.5}.} If $\gamma_k=1/L$ and $m_k=m$ for $k=1,\ldots,{N_{in}}$, \eqref{proj-gra-rand-bound} implies that
\bee
\E[\|\tilde{g}_R^r\|^2]  \le \frac{D_{\Phi_h} + {N_{in}}\sigma^2/(L m)}{{N_{in}}/(2L)} = \frac{2L D_{\Phi_h}}{{N_{in}}} + \frac{2\sigma^2}{m}.
\eee
Using Lemma \ref{cor2.1} with $g_1=G_k$ and $g_2=\nabla f(x_k)$, we have
\begin{align*}
\E[\|\tilde{g}_R\|^2] & \le  2 \E[\|\tilde{g}_R^r\|^2] + 2\E[\|\tilde{g}_R^r - \tilde{g}_R\|^2]  \le  \frac{4L D_{\Phi_h} }{{N_{in}}} + \frac{4\sigma^2}{m} + 2\E[\|G_R - \nabla f(x_R)\|^2]  \le   \frac{4L D_{\Phi_h} }{{N_{in}}} + \frac{6\sigma^2}{m},
\end{align*}
Note that the number of iterations of Algorithm \ref{alg-uncons} is at most ${N_{in}}=\lceil \bar{N}/m\rceil$. Obviously, ${N_{in}}\ge \bar{N}/(2m)$. Then following from \eqref{batch-size-m} we have that
\begin{align}
\E[\|\tilde{g}_R\|^2] & \le \frac{4L D_{\Phi_h}}{{N_{in}}} + \frac{6\sigma^2}{m}  \le  \frac{8L D_{\Phi_h}}{\bar{N}}m + \frac{6\sigma^2}{m} \label{exp-1} \\
& \le  \frac{8L D_{\Phi_h} }{\bar{N}}\left(1+\frac{\sigma}{L}\sqrt{\frac{\bar{N}}{\tilde{D}}}\right) + 6 \max\left\{\frac{ \sigma^2}{\bar{N}}, \frac{ \sigma L\sqrt{\tilde{D}}}{\sqrt{\bar{N}}}\right\}.  \label{exp-2}
\end{align}
From \eqref{bound-N} we have
\begin{align}
\sqrt{\bar{N}} & \ge \frac{\sqrt{(D_{\Phi_h}C_2 + LC_3)^2 + 32LD_{\Phi_h}\epsilon}}{\epsilon} \notag \\
 &\ge \frac{\sqrt{(D_{\Phi_h}C_2 + LC_3)^2 + 32LD_{\Phi_h}\epsilon} + (D_{\Phi_h}C_2 + L C_3)}{2\epsilon}  . \label{barN-bound}
\end{align}
\eqref{bound-N} also suggests that $\sigma^2/\bar{N}\le \sigma L \sqrt{\tilde{D}}/\sqrt{\bar{N}}$, which indicates from \eqref{exp-2} that
\begin{equation}\label{bound-g}
\E[\|\tilde{g}_R\|^2] \le \frac{8L D_{\Phi_h}}{\bar{N}} + \frac{8\sigma D_{\Phi_h}}{\sqrt{\bar{N}\tilde{D}}} + \frac{6L\sigma}{\sqrt{\bar{N}}}\sqrt{\tilde{D}} = \frac{8L D_{\Phi_h}}{\bar{N}} + \frac{D_{\Phi_h}C_2 + LC_3}{\sqrt{\bar{N}}} \le \epsilon,
\end{equation}
where the last inequality follows from \eqref{barN-bound}.
Note that \eqref{bound-g} together with \eqref{exp-1} implies that
\bee 
4L D_{\Phi_h}/{N_{in}} + 6\sigma^2/m\le \epsilon,
\eee
which according to \eqref{exp-1} shows that $\E[\|\tilde{g}_R^r\|^2]\le\epsilon$.
\qed

The following is the proof of Theorem \ref{thm4.1}.

{\it Proof of Theorem \ref{thm4.1}.}  It follows from part a) of Lemma \ref{lem4.1} that $f_\mu\in\mathcal{C}_{L_\mu}^{1,1}$ {with $L_\mu\le L_g$}. By {\bf AS.5}, \eqref{G}, \eqref{ine-3} and \eqref{G_mu} we obtain
\begin{align*}\label{G-mu-f}
\E_{v_k,\xi_k}[\|G_\mu(x_k,\xi_k,v_k) - \nabla f_\mu(x_k)\|^2] &\le \E_{v_k,\xi_k}[\|G_\mu(x_k,\xi_k,v_k)\|^2] \\
&\le  \E_{\xi_k} \left[2(n+4)\|G(x_k,\xi_k)\|^2 + \frac{\mu^2}{2} L_g^2(n+6)^3 \right] \\
& =  2(n+4) \E_{\xi_k}[\|G(x_k,\xi_k)\|^2] + \frac{\mu^2}{2} L_g^2(n+6)^3\\
& \le  2(n+4)(\|\nabla f(x_k)\|^2 + \sigma^2) + 2\mu^2L_g^2(n+4)^3 \le \tilde{\sigma}^2,
\end{align*}
where the last inequality follows from that {\bf AS.4} holds for $G(x_k,\xi_k)$. 
Similar to \eqref{delta-bound}, we can show that
\be \label{G-mu-f-zero}
\E[\|G_{\mu,k} - \nabla f_\mu(x_k)\|^2] \le  \frac{\tilde{\sigma}^2}{m_k}
\ee
according to the definition of $G_{\mu,k}$ in \eqref{G_mu_k}.

Denote
$ \Phi_{\mu,h}(x) := f_\mu(x) + h(c(x))$ and $\Phi_{\mu,h}^* = \min_{x\in\R^n} \Phi_{\mu,h}(x)$.
{\bf AS.3} {together with the continuity of $\Phi_{\mu,h}$} indicates that $\Phi_{\mu,h}^*$ is well-defined. So there exists $\hat{x}\in\R^n$ such that $\Phi_{\mu,h}^*=\Phi_{\mu,h}(\hat{x})$. By noting that $\Phi_{\mu,h}(x)-\Phi_h(x)=f_\mu(x)-f(x)$, we have from \eqref{ine-1} that
\begin{align*}
\Phi_{\mu,h}(x_1) - \Phi_{\mu,h}^* & = \Phi_{\mu,h}(x_1) - \Phi_{\mu,h}(\hat{x})\\
& =  \Phi_h(x_1) - \Phi_{h}(\hat{x}) +  \Phi_{\mu,h}(x_1) - \Phi_h(x_1)  - (\Phi_{\mu,h}(\hat{x})  - \Phi_{h}(\hat{x})) \\
& \le \Phi_h(x_1) - \Phi_h^{low} + |\Phi_{\mu,h}(x_1) - \Phi_h(x_1)| + | \Phi_{\mu,h}(\hat{x})-\Phi_h(\hat{x})| \\
& \le  \Phi_h(x_1) - \Phi_h^{low} + \mu^2 L_gn  \\
& = D_{\Phi_h} + \mu^2 L_gn.
\end{align*}
Therefore, by replacing $f$ with $f_\mu$ and $G_k$ with $G_{\mu,k}$ in Theorem \ref{thm2.2} we obtain
\begin{align*}
\E[\|\tilde{g}_{\mu,R}^r\|^2]   \le  \frac{\Phi_{\mu,h}(x_1) - \Phi_{\mu,h}^* + \tilde{\sigma}^2 \sum_{k=1}^N(\gamma_k/m_k)}{\sum_{k=1}^N(\gamma_k - L\gamma_k^2/2)} \le \frac{D_{\Phi_h} + \mu^2L_gn + \tilde{\sigma}^2 \sum_{k=1}^{N_{in}}(\gamma_k/m_k)}{\sum_{k=1}^{N_{in}}(\gamma_k - L\gamma_k^2/2)},
\end{align*}
where the expectation is taken with respect to $R$, $\xi_{[{N_{in}}]}$ and $v_{[{N_{in}}]}$.
\qed

We now give the proof of Theorem \ref{cor4,4}.

{\it Proof of Theorem \ref{cor4,4}.}\quad
 It follows directly from \eqref{proj-gra-rand-bound-zero} with $\gamma_k=1/L$ and $m_k=m$ that
 \bee
 \E[\|\tilde{g}_{\mu,R}^r\|^2]  \le  \frac{2L D_{\Phi_h} + 2\mu^2 LL_g n}{{N_{in}}} + \frac{2\tilde{\sigma}^2}{m}.
 \eee
 Note that
\begin{align}
\E[\|\tilde{g}_R\|^2] & \le  2\E[\|\tilde{g}_{\mu,R} - \tilde{g}_R\|^2] + 2\E[\|\tilde{g}_{\mu,R}\|^2]   \le  2\E[\|\tilde{g}_{\mu,R} - \tilde{g}_R\|^2] + 4\E[\|\tilde{g}_{\mu,R}^r\|^2] + 4\E[\|\tilde{g}_{\mu,R}^r - \tilde{g}_{\mu,R}\|^2]. \label{exp-tile-g}
\end{align}
Firstly, definitions of $\tilde{g}_{k}$ and $\tilde{g}_{\mu,k}$ in \eqref{g-rand_g} and \eqref{g-rand_g-zero} and Lemma \ref{cor2.1} indicate that
\bee
\|\tilde{g}_{\mu,R} - \tilde{g}_{R}\|^2 \le  \|\nabla f_\mu(x_R) - \nabla f(x_R)\|^2,
\eee
which together with \eqref{ine-2} shows that
\bee 
\|\tilde{g}_{\mu,R} - \tilde{g}_R\|^2 \le \frac{1}{4}\mu^2 L_g^2(n+3)^3.
\eee
Secondly, the definition of $\tilde{g}_{\mu,k}^r$ in \eqref{g-rand_g-zero} implies that
\be \label{diff-gmur-gmu}
\E[\|\tilde{g}_{\mu,R}^r - \tilde{g}_{\mu,R}\|^2]  \le   \E[\|G_{\mu,R} - \nabla f_\mu(x_R)\|^2]  \le \frac{\tilde{\sigma}^2}{m},
\ee
where the second inequality is due to \eqref{G-mu-f-zero}. Therefore, \eqref{exp-tile-g}-\eqref{diff-gmur-gmu} yield
\be\label{E_R-zero}
\E[\|\tilde{g}_R\|^2]  \le  \frac{1}{2}\mu^2 L_g^2(n+3)^3 + \frac{8L D_{\Phi_h} + 8\mu^2 LL_g n}{{N_{in}}} + \frac{8\tilde{\sigma}^2}{m} + \frac{4\tilde{\sigma}^2}{m}.
\ee
Given the total number of $\SZO$-calls $\bar{N}$ in the whole algorithm and the number of $\SZO$-calls $m$ at each iteration, we know that the inner iteration number of Algorithm \ref{alg-uncons-zero} is at most ${N_{in}}=\lceil\bar{N}/m \rceil\ge\bar{N}/(2m)$. Then \eqref{mu} and \eqref{E_R-zero} imply that
\begin{align}
\E[\|\tilde{g}_R\|^2] & \le  \frac{1}{2}\mu^2 L_g^2(n+3)^3 + \frac{16LD_{\Phi_h} + 16\mu^2LL_gn}{\bar{N}}m + \frac{12\tilde{\sigma}^2}{m}  \notag \\
& \le  \frac{\tilde{D}_1}{2\bar{N}} L_g^2(n+3)^3 + \frac{16LD_{\Phi_h}}{\bar{N}}m + \frac{ 16LL_gn}{\bar{N}}\cdot\frac{\tilde{D}_1}{\bar{N}}m + \frac{24(n+4)(\kappa_g^2+\sigma^2)}{m} + \frac{24(n+4)^3}{m}\cdot\frac{L_g^2\tilde{D}_1}{\bar{N}} \notag \\
& \le  \frac{25L_g^2\tilde{D}_1(n+4)^3 + 16LL_g\tilde{D}_1n}{\bar{N}} + \frac{16LD_{\Phi_h} }{\bar{N}}m + \frac{24(n+4)(\kappa_g^2+\sigma^2)}{m}  \notag\\
& \le  \frac{28LL_g\tilde{D}_1(n+4)^3}{\bar{N}} + \frac{16LD_{\Phi_h} }{\bar{N}}m + \frac{24(n+4)(\kappa_g^2+\sigma^2)}{m},\label{right-exp-R}
\end{align}
where we have used the fact that $1\le m\le \bar{N}$.
The choice of $m$ in \eqref{m} also yields that
\begin{align*}
\E[\|\tilde{g}_R\|^2] & \le \frac{28LL_g\tilde{D}_1(n+4)^3 }{\bar{N}} + \frac{16LD_{\Phi_h} }{\bar{N}}\left(1+\frac{1}{L}\cdot\sqrt{\frac{\bar{N}}{\tilde{D}_2}}\right)  + 24(n+4)(\kappa_g^2+\sigma^2)\cdot\max\left\{\frac{1}{\bar{N}}, \frac{L\sqrt{\tilde{D}_2}}{\sqrt{\bar{N}}}\right\}\\
 &= \frac{28LL_g\tilde{D}_1(n+4)^3 + 16LD_{\Phi_h} }{\bar{N}} + \frac{16D_{\Phi_h}}{\sqrt{\bar{N}\tilde{D}_2}} + \frac{24L}{\sqrt{\bar{N}}}(n+4)(\kappa_g^2+\sigma^2)\cdot\max\left\{\frac{1}{L\sqrt{\bar{N}}},\sqrt{\tilde{D}_2}\right\}.
\end{align*}
Then similar to the proof in Theorem \ref{cor2.5}, according to \eqref{bar_N-zero} it is easy to check that \eqref{eps-exp} holds.
\qed


\end{document}